\documentclass{article}

\usepackage[export]{adjustbox}
\usepackage{algpseudocode}
\usepackage{algorithm}
\usepackage{array}
\usepackage{amsthm}
\usepackage{amsmath}
\usepackage{amssymb}
\usepackage{booktabs}
\usepackage{comment}
\usepackage{etoolbox}
\usepackage{float}
\usepackage[T1]{fontenc}
\usepackage[top=1in, left=1in, right=1in, bottom=1in]{geometry}
\usepackage{gensymb}
\usepackage{graphicx}
\usepackage{hyperref}
    \hypersetup{colorlinks,allcolors=blue}
    \usepackage{cleveref}
\usepackage[utf8]{inputenc}
\usepackage{lineno}
\usepackage{mathtools}
\usepackage{multicol}
\usepackage{placeins}
\usepackage{subcaption}

\usepackage{tikz}
    \usetikzlibrary{calc}
    \usetikzlibrary{external}
    \usetikzlibrary{quotes}
    \usetikzlibrary{chains}
    \tikzset{vertex/.style={black,fill,draw,minimum size=6pt,inner sep=0pt,circle,thin},bold/.style={black,line width=0.6mm},plain/.style={black,thin},bold edges/.style=bold,plain edges/.style=plain,label distance=1mm,text node/.style={rectangle,fill=none,draw=none},every label/.style=text node,caption node/.style={text node,font=\Large}}

\usepackage[dvipsnames]{xcolor}
\usepackage{paralist}
\newtheorem{theorem}{Theorem}[section]
\newtheorem{lemma}[theorem]{Lemma}
\newtheorem{corollary}[theorem]{Corollary}

\newtheorem{proposition}[theorem]{Proposition}

\newtheorem{question}[theorem]{Question}
\newtheorem{example}[theorem]{Example}

\newtheorem{remark}{Remark}

\newtheorem{definition}[theorem]{Definition}

\newcommand{\abs}[1]{\left| #1 \right|}
\newcommand{\ceil}[1]{\left\lceil #1 \right\rceil}
\newcommand{\ar}[1]{\overset{\rightarrow}{#1}}
\newcommand{\AR}[1]{\overrightarrow{#1}}

\newcommand{\ra}[1]{{\color{olive}#1}}

\title{Questions of minimality for DDMOGs and cycle ornaments}
\author{Roza Aceska, McKailyn Lort, Allison Ripperger and Zach Wilcher}
\date{today}

\begin{document}

\maketitle

\begin{abstract}
Difference distance magic oriented graphs (DDMOGs) are oriented graphs whose vertices admit a bijective labeling by consecutive positive integers such that every vertex has weight zero, where the weight is the sum of the labels of its in-neighbors minus the sum of the labels of its out-neighbors. We investigate minimality and sparsity questions for DDMOGs and their ornaments. First, we develop new constructions of cycle-based ornaments that reduce the number of cycles required for a given number of attachment nodes and identify minimal cycle lengths in several cases. 
Second, we determine the minimum number of arcs in a nontrivial DDMOG of order $n$, combining computational enumeration for small orders with theoretical lower bounds and constructions. Langford sequences provide a key tool in these constructions, enabling us to partition consecutive labels into difference triples and thereby construct sparse DDMOGs. Finally, we use related constructions to obtain generalized ornaments with no internal arcs. These results advance the study of minimal constructions in difference distance magic labeling and provide new methods for building sparse DDMOGs.


\end{abstract}
\section{Introduction}


Graph labeling is an assignment of values to the components of a graph (its vertices, edges, etc.) such that a certain property holds \cite{ BB}. Various types of labelings exist and Gallian's survey is one source to keep updated on what has been done in the field \cite{Gal}. Since 2003, distance magic labelings have been extensively  studied on undirected graphs,  \cite{undi}. In such a labeling, the vertices are labeled with unique values so that the sum of labels on the neighbors of any vertex is the same across all vertices.   Difference distance magic (DDM) labelings  on oriented graphs were initially introduced in  \cite{FM},  and further explored in  \cite{AAetc}, where the authors offered a technique that creates a large(r) graph with  DDM labeling via a weighted sum of a smaller DDMOG with a well-selected oriented graph. In \cite{ALR} this graph was named an ornament; the authors of \cite{ALR} presented a collection of ornaments (each formed by four cycles with $L\geq 3$ vertices). The ornaments were categorized by the specifics of the utilized weighted sum.


\subsection{Definitions and known results}
Every graph in this work will be finite and simple (no loops and no multiple edges between adjacent vertices). Given a graph \(G\), we will primarily be labeling the vertices of an \textit{orientation} of \(G\). A directed graph \(\ar{G}\) is an orientation of \(G\) if and only if \(\ar{G}\) has the same vertex set as \(G\), and if every arc in \(\ar{G}\) is treated as an undirected edge, then the result is equal to \(G\). If \(\ar{G}\) is an orientation of some graph, then we say \(\ar{G}\) is an oriented graph.

\begin{definition}
\label{definition:vertex-weight}
Suppose \(\ar{G}\) is an oriented graph and \(g: V(\ar{G}) \rightarrow \mathbb{Z}^+\) is a vertex labeling of \(\ar{G}\). The \textit{weight} of a vertex \(v\) in \(\ar{G}\) with respect to \(g\) is
\(\displaystyle
wt_g(v) = \,\sum_{\mathclap{w \in N^{+}(v)}} \,g(w) - \,\sum_{\mathclap{w \in N^{-}(v)}} \,g(w).
\)
Here \(N^{+}(v)\) and \(N^{-}(v)\) respectively denote the incoming and outgoing neighborhoods of \(v\) in \(\ar{G}\). If the labeling is clear from context we will write the weight of \(v\) as \(wt(v)\) instead of \(wt_g(v)\).
\end{definition}

Definition \ref{definition:vertex-weight} specifies that the vertex labeling of \(\ar{G}\) is a map with values in \(\mathbb{Z}^+\). 
We restrict our focus to bijective labelings using the positive integers but note that others \cite{MR3703024, MR4926298} have studied distance magic labelings of oriented graphs which use finite Abelian groups.
We now give the definition of a DDMOG used by \cite{FM} in Definition \ref{definition:DDMOG} below.
\begin{definition}
\label{definition:DDMOG}
\cite{FM}
Suppose \(G\) is a graph with \(n\) vertices. We say \(G\) is \textit{difference distance magic orientable} (DDMO) if and only if there exists a bijective labeling \(g: V(G) \rightarrow \{1,2,\cdots,n\}\) and an orientation \(\ar{G}\) of \(G\) such that for some (magic) constant \(k\) every vertex \(v\) in \(\ar{G}\) has weight \(wt_g(v) = k\).
If \(\ar{G}\) is an orientation of \(G\) and \(g\) is a vertex labeling of \(G\) that makes \(G\) DDMO, then we say \(\ar{G}\) is a difference distance magic oriented graph (DDMOG) and \(g\) is a difference distance magic (DDM) labeling of \(\ar{G}\).
\end{definition}

Every graph without any edges is trivially DDMO {with magic constant 0}. In general, the magic constant of a DDMOG is always 0. If a DDMOG has any isolated vertices, then we say it is a \textit{trivial} DDMOG. Although the work in \cite{FM}  focused primarily on connected graphs, many of their results apply more generally. 


\begin{lemma}
\label{lemma:minimal-incident-arcs}
\cite{FM} Every vertex \(v\) in a nontrivial DDMOG is incident to at least 3 arcs, not all directed into or away from \(v\).
\end{lemma}

It follows from Lemma \ref{lemma:minimal-incident-arcs} that there are no nontrivial DDMOGs with less than 4 vertices. It is also shown in \cite{FM} that \(K_n\) is not DDMO for any \(n \ge 2\). Since any nontrivial DDMOG with 4 vertices must be isomorphic to \(K_4\), there are no nontrivial DDMOGs with 4 vertices. Hence, every nontrivial DDMOG has at least 5 vertices. Figure \ref{fig:a}, 
illustrates an example of a DDMOG with 5 vertices.


Lemma \ref{lemma:minimal-incident-arcs} also guarantees
that cycle graphs are never DDMO. However, this does not mean cycle graphs cannot be used as parts to construct DDMOGs. In \cite{AAetc} cycle graphs are attached to another graph using a new graph sum they call the \textit{weighted sum}.

\begin{definition}
\cite{AAetc}
\label{definition:weighted-sum}
For an oriented graph $\ar{G}$ with labeling $g: V(\ar{G}) \rightarrow \mathbb{Z}^+$, we write $V_g^i(\ar{G}) = \{v \in V(\ar{G}) : \text{wt}_f(v) = i\}$.
Let $k \in \mathbb Z^+$ and let $\ar{G}$ be an oriented graph with vertex set $V(\ar{G}) = \{v_1,\dots, v_n\}$ and associated labeling function $g: V(\ar{G}) \rightarrow \{1,2,\dots, n\}$, where $g(v_i)=i$. 

If $\ar{H}$ is an oriented graph with labeling $h: V(H) \rightarrow \mathbb{Z}^+$, where  
$\{ |wt_h(u)| : u \in V(\ar{H}) \} \subseteq \{0\} \cup \{ g(v)+k : v \in V(\ar{G})\}$ for some $k \in \mathbb Z$,  then the weighted sum of $\ar{G}$ and $\ar{H}$ shifted by $k$ denoted $\ar{G} \oplus_{wt_h}^k \ar{H}$ is the oriented graph with vertex set $V(\AR{G} \oplus_{wt_h}^k \AR{H})=V(\AR{G}) \cup V(\AR{H})$ and $E(\AR{G} \oplus_{wt_h}^k \AR{H})=E(\AR{G}) \cup E(\AR{H}) \cup (\bigcup\limits_{i=1}^n (E^i \cup E^{-i}))$, where $E^i=\{v_i u:u \in V_h^{-i-k}(\AR{H})\}$ and $E^{-i}=\{uv_i:u \in V_h^{i+k}(\AR{H})\}$; here $V_h^i(\AR{G})=\{v \in V(\AR{G}) \, : \, wt_h(v)=i\}$ and $V_h^{-i}(\AR{G})=\{v \in V(G) \, : \,  wt_h(v)=-i\}$.
\end{definition}

Figure \ref{fig:a} illustrates how a DDMOG can be viewed as the weighted sum of \(\ar{K_1}\) and \(\ar{C_4}\).
The weighted sum is offered as a way to construct new DDMOGs from known DDMOGs with \(n\) vertices, when paired with suitably selected graphs.  For different values of integers \(L\) and \(s\), the weighted sum is used in \cite{ALR} to \textit{attach} \(4\ar{C_L}\) to exactly \(s\) vertices in an arbitrary DDMOG, and thus create a new DDMOG. 

%
\begin{definition}
    \label{definition:ornament}
    \cite{ALR}
    If \(H\) and \(G\) are vertex-disjoint graphs, then we say \(H\) is attached to \(G\) if an edge joins them, and call the vertices in \(G\) that are adjacent to \(H\) the attachment nodes of \(H\).

    Let \(\ar{G}\) be a DDMOG with a DDM labeling \(f_0: V(\ar{G}) \rightarrow \{1,\cdots,|{V(\ar{G})}|\}\).
    An oriented graph \(\ar{H}\) is an \(s\)-node ornament attachable to \(\ar{G}\) if and only if there exists a bijective labeling \(h: V(\ar{H}) \rightarrow \{n+1,\cdots, n+|{V(\ar{H})}|\}\) 
    such that the function \(f: V(\ar{G})\cup V(\ar{H}) \rightarrow \{1,\cdots,|{V(\ar{G})}| + |{V(\ar{H})}|\}\) given by \begin{equation}\label{deffonG0H} f(v)=\begin{cases}
			h(v), & \text{if $v \in V(H)$,}\\
            f_0(v), & \text{if $v \in V(G_0)$.}
		 \end{cases}\end{equation}
    is a DDM labeling of the graph \(\ar{G} \oplus_{wt_h} \ar{H}\),
    and \(\ar{H}\) has exactly \(s\) attachments nodes. That is, exactly \(s\) vertices in \(\ar{G}\) are adjacent to vertices in \(\ar{H}\).
\end{definition}

Suppose \(\ar{H}\) is an oriented graph guaranteed to be an \(s\)-node ornament attachable to a DDMOG \(\ar{G_0}\) with \(n_0\) vertices by a labeling \(h\) of \(\ar{H}\) and a collection of arcs \(E_0\) attaching \(\ar{H}\) to \(\ar{G_0}\).
We emphasize that the only property required of \(\ar{G_0}\) for \(\ar{H}\) to attach to \(\ar{G_0}\) is the presence of the \(s\) vertex labels that \(\ar{H}\) attaches to.
So, any DDMOG \(\ar{G_1}\) with \(n_1 \ge n_0\) vertices can be used in place of \(\ar{G_0}\) in Definition \ref{definition:ornament} by using a set of arcs that attach labels in \(\ar{H}\) to labels in \(\ar{G_1}\) in the same fashion that the arcs in \(E_0\) attach \(\ar{H}\) to \(\ar{G_0}\).


The following theorem summarizes various ways \(4\ar{C_L}\) can be labeled to make it an \(s\)-node ornament.

\begin{theorem} 
\label{theorem:4CL-ornaments}
\cite{ALR}
    Let $\ar{G_0}$ be a DDMOG with $n$ vertices, and let  $s \in \{2,3, 4,5,6\}$. Let $L$ be an integer such that  $a_s\leq L\le b_s$, where    $a_s$, $b_s$ take values as listed in Table~\ref{s-table}, with the constraint that $L$ is even when $s=3$.
    If $H=4 C_L$ is a graph of  four $L-$cycles, then there exists  an orientation on $ {H} $, paired with  a labeling function $h$, such that 

     \begin{itemize}
         \item[a.] the oriented graph $ \ar{H}$ is a $s-$nodes ornament, where the labels of the nodes are as listed in Table~\ref{s-table}, and 
      \item[b.]  $\displaystyle  \ar{G}:=\ar{G_0}\oplus_{\text{wt}_h} \ar{H}$  is a DDMOG.
      \end{itemize}

      \begin{table}[h!]
 \centering
\begin{tabular}{|c|c|c|c|}
 \hline
 $s$ & $a_s$ & $b_s$ & Nodes labels \\
 \hline
 2 & 3 & $n+2$ & 2, $L-2$ \\
 \hline
 3 & 6 & $n+5$ & 4, $L-5$, $L-3$ \\
 \hline
 4 & 6 & $n+2$ & 1, 2, 3, $L-2$ \\
 \hline
 5 & 7 & $n+3$ & 1, 2, 3, $L-3$, $L-2$ \\
 \hline
 6 & 8 & $n+3$ & 1, 2, 3, 4, $L-3$, $L-2$ \\
 \hline
\end{tabular}
\caption{Node labels and bounds for $L$ in 
Theorem~\ref{theorem:4CL-ornaments}}
\label{s-table}
\end{table}
  
\end{theorem}
We use Theorem \ref{thmWgt1} below frequently in our constructions.

\begin{theorem}
\label{thmWgt1} \cite{AAetc} Let $\AR{G}$ be a DDMOG on $n$ vertices with DDM labeling $g$ and let $\AR{H}$ be an oriented graph on $m$ vertices with a bijective labeling $h: V(\AR{H}) \rightarrow \{n+1,\dots, n+m\}$ such that $k=\displaystyle \max_{v \in V(\AR{H})} \{|wt_h(v)|\} \leq n$ and for every $0 \leq i \leq k$, $\displaystyle \sum_{v \in V_h^i(\AR{H})} h(v)=\sum_{v \in V_h^{-i}(\AR{H})} h(v)$.  Then, $\AR{G} \oplus_{wt_h} \AR{H}$ is a DDMOG on $n+m$ vertices. 
\end{theorem}

In \cite{FM} several other DDMOG classes were considered. Some relevant facts include necessary and sufficient criteria for \(K_{p,q}\) to be DDMO (see Theorem \ref{theorem:K_p,q-is-DDMO}) as well as a complete classification of all nontrivial DDMOGs with \(n = 5\) and \(n = 6\) vertices.

\begin{theorem}
\label{theorem:K_p,q-is-DDMO}
\cite{FM} \(K_{p,q}\) is DDMO if and only if \(p,q \ge 3\) and \((p + q)(p + q + 1) \equiv 0 \pmod{4}\).
\end{theorem}

To verify their classification of 5 and 6 vertex DDMOGs, the authors in \cite{FM} used a program to iterate through every nontrivial DDMOGs with \(n\) vertices. This program successfully found the 6 nontrivial DDMOGs with 5 vertices and the 22 nontrivial DDMOGs with 6 vertices. 
Their program also found that there are 296 nontrivial DDMOGs with 7 vertices and 8240 nontrivial DDMOGs with 8 vertices. Interested in finding nontrivial DDMOGs with the fewest number of arcs possible, we also used a program to search for every possible nontrivial DDMOG with \(n\) vertices (see Algorithm \ref{alg:DDMOGIterator}). Table \ref{table:DDMOG-counts} contains the results of our program's computations. 

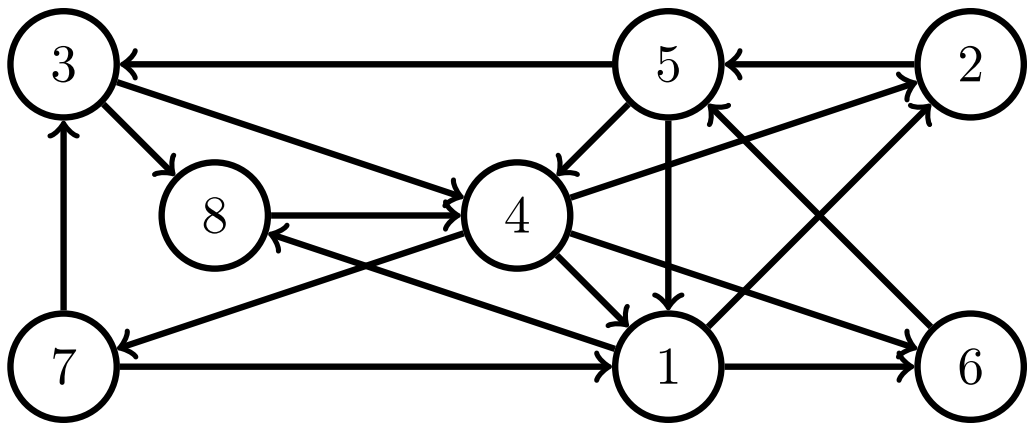
\begin{figure}[h!]
\centering
\begin{tikzpicture}[scale=2, rotate=90]
    \tikzset{
        labeled_vertex/.style = {circle, draw=black, fill=white, inner sep = 1.5pt, text=black, minimum size=20pt, very thick},
        arc/.style = {->, very thick}
    };

    \node [labeled_vertex] (6) at (0, 0) {6};
    \node [labeled_vertex] (2) at (1, 0) {2};
    \node [labeled_vertex] (1) at (0, 1) {1};
    \node [labeled_vertex] (5) at (1, 1) {5};
    \node [labeled_vertex] (4) at (0.5, 1.5) {4};
    \node [labeled_vertex] (8) at (0.5, 2.5) {8};
    \node [labeled_vertex] (7) at (0, 3) {7};
    \node [labeled_vertex] (3) at (1, 3) {3};

    \draw[arc] (1) to (2);
    \draw[arc] (4) to (1);
    \draw[arc] (5) to (1);
    \draw[arc] (1) to (6);
    \draw[arc] (7) to (1);
    \draw[arc] (1) to (8);

    \draw[arc] (4) to (2);
    \draw[arc] (2) to (5);

    \draw[arc] (3) to (4);
    \draw[arc] (5) to (3);
    \draw[arc] (7) to (3);
    \draw[arc] (3) to (8);

    \draw[arc] (5) to (4);
    \draw[arc] (4) to (6);
    \draw[arc] (4) to (7);
    \draw[arc] (8) to (4);

    \draw[arc] (6) to (5);
    
\end{tikzpicture}
\caption{An order 8 DDMOG not counted in \cite{FM}}
\label{fig:order-8-non-muell-DDMOG}
\end{figure}

Note that our count for the number of DDMOGs with \(n = 8\) vertices is slightly larger. This discrepancy is due to certain DDMOGs being excluded from their count. In Figure \ref{fig:order-8-non-muell-DDMOG} we illustrate an example of this. Notice that the weight of the vertex labeled with 4 in the figure is
\(-1 -2 + 3 + 5 - 7 + 8\).
Unfortunately once the algorithm used in \cite{FM}  determines that \(-1 - 2 + 3 = 0\), it does not consider the possibility that other labels can be included along with -1, -2, and 3 in 4's weight equation and the total still be 0.

\subsection{Questions of minimality}

\begin{question}
    Given integers \(L \ge 3\) and \(s \ge 1\), what is the minimal \(p\) such that \(p\ar{C_L}\) is an \(s\)-node ornament?
\end{question}

In Section \ref{section:minimal-p},
Theorem \ref{thm:2nodes1cycle} shows that \(\ar{C_L}\) can function as a \(2\)-node ornament by using weighted sums, but the construction requires that \(L\) is a multiple of 4. Theorem \ref{thm:3nodes2cycles} shows that \(2\ar{C_L}\) can function as a \(3\)-node ornament when \(L\) is a multiple of 2; see Figure~\ref{fig:Two6Cycles} for a detailed construction when $L=6$. In Section \ref{section:minimal-L}, we flip this question and try to minimize \(L\) instead. 

\begin{question}
    Given positive integers \(p\) and \(s\), what is the minimal value of \(L\) such that \(p\ar{C_L}\) is an \(s\)-node ornament?
\end{question}

Since all cycles have at least 3 vertices, Theorem \ref{theorem:4CL-ornaments} minimizes \(L\) when \(s = 2\) with a labeling of \(4\ar{C_3}\). However, if \(s = 1\) then \(L\) must be exactly 4 (see Proposition \ref{prop:c4is1nodeonly}). 
Figure \ref{fig:a} shows that \(\ar{C_4}\) is a 1-node ornament, and in the discussion leading up to Corollary \ref{corollary:jC_4-is-1-node-ornament} we show that \(p\ar{C_4}\) is a 1-node ornament for any \(p > 1\).

Definition \ref{definition:ornament} does not require that an ornament is of the form \(p\ar{C_L}\). As suggested in \cite{ALR}, if this restriction is removed, then another question of minimality arises.

\begin{question}
\cite{ALR}
Given positive integers \(s\) and \(k\), what is the fewest number of arcs an \(s\)-node ornament with \(k\) vertices can have?
\end{question}
In Section \ref{section:sparse-ddmog-search} we provide a partial answer to this question in
Proposition \ref{proposition:3mK_1-ornament} which shows that \(3mK_1\) can function as an ornament for sufficiently large \(m\). The strategy employed in Proposition \ref{proposition:3mK_1-ornament} arose exploring another question posed in \cite{ALR}, namely:
\begin{question}
\label{question:minimal-sparsity}
\cite{ALR}
    Given a positive integer \(n\), what is the fewest number of arcs a nontrivial DDMOG can have?
\end{question}
Lemma \ref{lemma:minimal-incident-arcs} states there needs to be at least 3 arcs incident to each vertex, not all directed into or away from it. 
Then, the degree sum formula guarantees that every nontrivial DDMOG must have at least \(\ceil{3n/2}\) arcs. 
Theorem \ref{theorem:4CL-ornaments} shows that \(4C_3 \oplus 2K_1\) is DDMO with 24 arcs only three more guaranteed by this lower bound.
It turns out \(4C_3 \oplus 2K_1\) is not ``minimally sparse.'' Figure \ref{fig:order-14-sparse-connected-DDMOG} and Figure \ref{fig:order-14-sparse-disconnected-DDMOG}.
both illustrate DDMOGs with 14 vertices and \(21 = 3\cdot 14 / 2 + 1\) arcs which Lemma \ref{lemma:minimum-edges} guarantees is minimal. We define precisely what we mean by minimally sparse in Definition \ref{definition:minimal-sparsity} below.

\begin{definition}
    \label{definition:minimal-sparsity}
    A nontrivial DDMOG \(\ar{G}\) with \(n\) vertices has minimal sparsity if and only if \(\vert E(\ar{G}) \vert \le \vert E(\ar{H})\vert \) for every nontrivial DDMOG \(\ar{H}\) with \(n\) vertices. If there is a nontrivial DDMOG with \(n\) vertices, \(MS(n)\) denotes the number of arcs in a minimal sparsity nontrival DDMOG with \(n\) vertices.
\end{definition}
With this notation, Question \ref{question:minimal-sparsity} asks ``What is the value of \(MS(n)\) for \(n \ge 1\)?''
However, this question only makes sense if \(n \ge 5\), 
since there are no nontrivial DDMOGS with \(n\) vertices when \(n < 5\). 

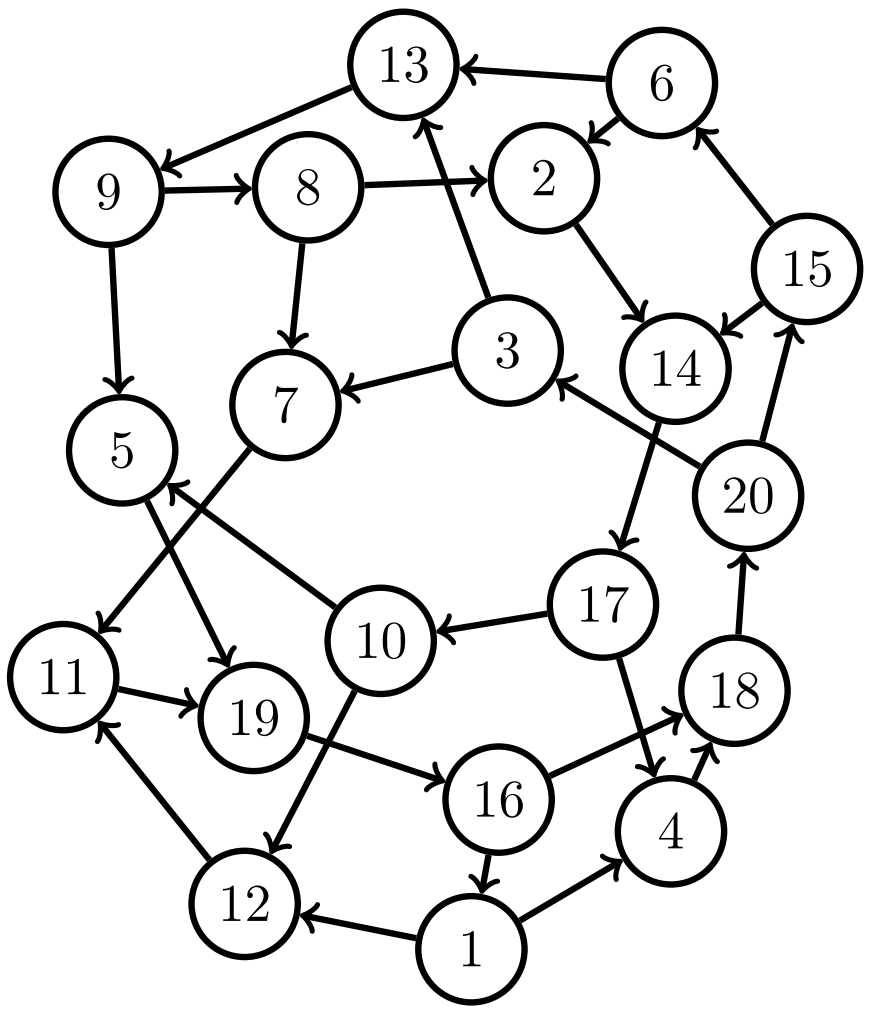
\begin{figure}[h!]
\centering
\begin{tikzpicture}[scale=3]
    \tikzset{
        labeled_vertex/.style = {circle, draw=black, fill=white, inner sep = 1.5pt, text=black, minimum size=20pt, very thick},
        arc/.style = {->, very thick}
    };
    \node[labeled_vertex] (1) at (0.00, -1.00) {1};
    \node[labeled_vertex] (2) at (0.16, 0.70) {2};
    \node[labeled_vertex] (3) at (0.08, 0.32) {3};
    \node[labeled_vertex] (4) at (0.44, -0.74) {4};
    \node[labeled_vertex] (5) at (-0.77, 0.1) {5};
    \node[labeled_vertex] (6) at (0.42, 0.91) {6};
    \node[labeled_vertex] (7) at (-0.41, 0.20) {7};
    \node[labeled_vertex] (8) at (-0.36, 0.68) {8};
    \node[labeled_vertex] (9) at (-0.8, 0.67) {9};
    \node[labeled_vertex] (10) at (-0.20, -0.32) {10};
    \node[labeled_vertex] (11) at (-0.9, -.4) {11};
    \node[labeled_vertex] (12) at (-0.5, -0.9) {12};
    \node[labeled_vertex] (13) at (-.15, 0.95) {13};
    \node[labeled_vertex] (14) at (0.45, 0.28) {14};
    \node[labeled_vertex] (15) at (0.74, 0.50) {15};
    \node[labeled_vertex] (16) at (0.06, -0.67) {16};
    \node[labeled_vertex] (17) at (0.29, -0.24) {17};
    \node[labeled_vertex] (18) at (0.58, -0.43) {18};
    \node[labeled_vertex] (19) at (-0.48, -0.49) {19};
    \node[labeled_vertex] (20) at (0.61, 0.00) {20};
    \draw[arc] (1) to (4);
    \draw[arc] (1) to (12);
    \draw[arc] (2) to (14);
    \draw[arc] (3) to (7);
    \draw[arc] (3) to (13);
    \draw[arc] (4) to (18);
    \draw[arc] (5) to (19);
    \draw[arc] (6) to (2);
    \draw[arc] (6) to (13);
    \draw[arc] (7) to (11);
    \draw[arc] (8) to (2);
    \draw[arc] (8) to (7);
    \draw[arc] (9) to (5);
    \draw[arc] (9) to (8);
    \draw[arc] (10) to (5);
    \draw[arc] (10) to (12);
    \draw[arc] (11) to (19);
    \draw[arc] (12) to (11);
    \draw[arc] (13) to (9);
    \draw[arc] (14) to (17);
    \draw[arc] (15) to (6);
    \draw[arc] (15) to (14);
    \draw[arc] (16) to (1);
    \draw[arc] (16) to (18);
    \draw[arc] (17) to (4);
    \draw[arc] (17) to (10);
    \draw[arc] (18) to (20);
    \draw[arc] (19) to (16);
    \draw[arc] (20) to (3);
    \draw[arc] (20) to (15);
\end{tikzpicture}
\caption{A 20 vertex DDMOG with 30 edges. Since \(30 = \ceil{\frac{3\cdot20}{2}}\), Lemma \ref{lemma:minimum-edges} guarantees this DDMOG has minimal sparsity.}
\label{fig:order-20-sparse-connected-DDMOG}
\end{figure}

\begin{theorem}
\label{theorem:minimal-sparsity}

The values of $MS(n)$ are  as listed in Table~\ref{table:minimal-sparsity-counts}. 
 
\begin{table}[h!]
 \centering
\begin{tabular}{ c | c c c c c c c}

 $n$ & 5 & 6 & 7 & 8 & 9 & $n \ge 10$, $n \equiv 2 \pmod{4}$ & $n \ge 10$, $n \not \equiv 2 \pmod{4}$\\
\hline
\textbf{$MS(n)$} & 8 & 11 & 12 & 13 & 15 & $3n/2+1$ & $\lceil 3n/2 \rceil$ \\

\end{tabular}
\caption{Values of $MS(n)$}
\label{table:minimal-sparsity-counts}
\end{table}

\end{theorem}
\begin{proof}
Table \ref{table:DDMOG-counts} summarizes the results of a brute force computation of \(MS(n)\) for \(n \le 10\).
Lemma \ref{lemma:minimum-edges} gives a lower bound for \(MS(n)\) which the examples in Figure \ref{fig:order-20-sparse-connected-DDMOG}, Proposition \ref{proposition:minimal-sparsity-families} and Theorem \ref{theorem:G2xK_3,3isDDMO} attain.
Specifically, Proposition \ref{proposition:minimal-sparsity-families} handles the case when \(n > 10\), \(n \neq 20\), and \(n \not \equiv 0 \pmod{12}\).
Figure \ref{fig:order-20-sparse-connected-DDMOG} illustrates a DDMOG with \(n = 20\) vertices that attains the lower bound in Lemma \ref{lemma:minimum-edges}.
Theorem \ref{theorem:G2xK_3,3isDDMO} where \(G\) is the graph with zero vertices handles the case when \(n \equiv 0 \pmod{12}\).
\end{proof}

Proposition \ref{proposition:minimal-sparsity-families} gives examples of DDMOGs whose size attains the value of \(MS(n)\) in Theorem \ref{theorem:minimal-sparsity} for most values of \(n \ge 10\). Although, the constructed graphs are generally disconnected. We were able to find connected examples that attain the value of \(MS(n)\) when \(13 \le n \le 23\), but not when \(n = 11, 12\), or when \(n > 23\).
Additionally, we were unable to find an example of a disconnected DDMOG with \(20\) vertices whose size attains the value of \(MS(20)\).

\begin{question}
Given any positive integer \(n\),
how many arcs must a connected minimal sparsity DDMOG have?
Moreover, for which \(n\) does this number of arcs differ from \(MS(n)\)?
\end{question}

We emphasize that Proposition \ref{proposition:minimal-sparsity-families} uses Theorem \ref{theorem:G2xK_3,3isDDMO}.
Given some DDMOG \(G\), this Theorem specifies a way to provide a DDM labeling of a graph of the form \(G\sqcup 2xK_{3,3}\) that preserves the labeling of \(G\). A particular orientation of \(K_{3,3}\) allows it to easily receive a DDM labeling using two \textit{difference triples} (a tuple of the form \((a,b,c)\) where \(a + b = c\)). By partitioning an array of the form \((n + 1, n + 2, \cdots, n + 12x)\) into difference triples, we can construct a DDM labeling of \(G\sqcup 2xK_{3,3}\) by extending \(G\)'s DDM labeling with difference triples as the labels of the \(2xK_{3,3}\) components. An order \(4x\) \textit{Langford pairing} with defect \(n+1\) determines a partition of the array \((n+1, n+2, \cdots, n+ 12x)\) into difference triples (see Lemma \ref{lemma:langford_triples} for the specific construction).
We introduce Langford pairings now.

\begin{definition}
    Let \(n\) and \(d\) be positive integers.
    An order \(n\) Langford sequence with defect \(d\)
    is an array \((l_1, l_2, \cdots, l_{2n})\)
    such that:
    \begin{enumerate}
        \item 
    For every \(t \in \{d, d+1, d+2, \cdots d+n-1\}\) there exists exactly two elements \(l_i\) and \(l_j\) where \(l_i = l_j = t\).
    
        \item If \(l_i = l_j = t\) and \(i < j\), then \(j - i = t\).
    \end{enumerate}
Often, it is more convenient to refer to the indices of a Langford sequence as a collection of ordered pairs 
\[
\{(u_i, v_i) : 1 \le i \le n, v_i - u_i = i + d - 1\}
\]
where \(\bigcup_{i=1}^n \{u_i, v_i\} = \{1,2,\cdots,2n\}\). We call such a collection of pairs an order \(n\) Langford pairing with defect \(d\).
\end{definition}

This definition is adapted from Definition 53.6 in \cite{colbournHandbookCombinatorialDesigns2006}.
Note that an order \(n\) Langford sequence with defect \(1\)
is also commonly known as an order \(n\) Skolem sequence.

\begin{example}
(3, 4, 2, 3, 2, 4) is a Langford sequence with order 3 and defect 2 and its corresponding Langford pairing is \(\{(3,5), (1, 4), (2, 6)\}\).
\end{example}

\begin{example}
(1, 1, 3, 4, 5, 3, 2, 4, 2, 5) is an order 5 Skolem sequence and its corresponding pairing is \(\{(1,2), (7, 9), (3,6), (4,8), (5,10)\}\)
\end{example}

Langford pairings exist with arbitrarily large defects, and there are even constructive formulas for them.

\begin{lemma}
    \label{lemma:langford-sequence-existence}
    \cite{simpsonLangfordSequencesPerfect1983}
    An order \(n\) Langford sequence with defect \(d\)
    exists if and only if
    \(n \ge 2d - 1\), and if \(d\) is odd, then \(n \equiv 0,1\pmod{4}\), otherwise \(n\equiv 0,3\pmod{4}\).

\end{lemma}

\subsection{Outline}
The remainder of this paper is organizes as follows: Section~\ref{section:minimal-p} explores cycle-based ornaments,  with several results on minimal cycle counts (Subsection~\ref{3nodesExample}, Subsection~\ref{subs:onecycle}) and cycle lengths (Subsection~\ref{subs:minLornament}, Subsubsection~\ref{subsection6C4}) in $s-node$ ornaments, $1\leq s \leq 6$. Figure~\ref{fig:Two6Cycles} illustrates in detail one of these results. The main focus of Section~\ref{section:sparse-ddmog-search} is graph sparsity, where tight lower bounds on the minimum arc count $MS(n)$ for nontrivial DDMOGs  are presented. What is more, Langford pairings are employed to construct an infinite family of minimally sparse graphs.   Finally in Section~\ref{sec:final},  instead of a classical conclusion, we use related constructions to introduce generalized ornaments with no internal arcs, that is, zero arcs ornaments (of type  \(3mK_1\)).

\section{Atypical or Minimal cycle ornaments}
\label{section:minimal-p}
By Theorem~\ref{theorem:4CL-ornaments}, four copies are typically needed to create an $s-$nodes ornament. We note that the authors of \cite{ALR} were using a repetitive labeling pattern on the cycles involved; in this section we will occasionally abandon that repetitive style of labeling.  {We first present special cases of cycle ornaments, where fewer than four cycles are needed - in other words, we aim to minimize the value of $p$ in ornaments of type $ p\ar{C_L}$, given any $L\geq 3$.  We will also explore the minimal value of $L$ such that $p\ar{C_L}$ is an $s-$nodes ornament, for some fixed values of $p$ and $s$. For this purpose we will also explore variations in labeling when $L=4$.

 \subsection{ {A} $3-$nodes ornament, formed   by  two $L-$cycles}\label{3nodesExample}
      
Under certain constraints, two oriented $L-$cycles form a $3-$node ornament. An example  of this type of ornament is shown in 
Subfigure~\ref{fig:b}: the ornament consists of two oriented $6-$cycles ($2\ar{C_6}$), with vertices  labeled $6$ through $ 17$. In the subfigure,  the respective weights (induced by the labeling on the cycles, $h$) are marked in orange, red  and teal. 
 Since  $2\ar{ C_6}$ is a $3-$node ornament, we can employ it to create a new DDMOG, starting with a known DDMOG, $\ar{G_0}=\ar{W}_4$ (Subfigure~\ref{fig:a}), with related DDM labeling $f_0$. The  vertices of ${G_0}={W}_4$ are assigned with labels $1$ through $ 5$. 
%
\noindent  
By introducing  new edges to and from the vertices labeled with $1, 3, 4$, we have created a connected oriented  graph $\vec{G} : = \vec{W_4} \oplus_{\text{wt}_f} (\vec{H_1} \cup \vec{H_2})$, where $f$  is the labeling function  as depicted in  Subfigure~\ref{fig:c}, same as in equation~\eqref{deffonG0H} (with $\ar{G_0}=\vec{W_4}$, $\ar{H}=2 \vec{C_6}$).

We now verify  that $\ar{G} $ is a DDMOG: First, it is trivial to see that $f : V({G}) \rightarrow \{1,2, \hdots, 17\}$ is a bijection. The weight $\text{wt}_f$ is zero at every vertex  $u \in V(2\vec{C_6})$, due to the newly added edges.
As for the weight values at the vertices of $\vec{W_4}$, we already have $\text{wt}_{f_0}(v)=0$ for every $v \in V(W_4)$, so   $\text{wt}_f(2)=\text{wt}_{f_0}(2)=0$, $\text{wt}_f(5)=\text{wt}_{f_0}(5)=0$. At the three nodes (labeled  $1$, $3$, $4$) we have:  $\text{wt}_f(1) =\text{wt}_{f_0}(1) + 17-12 + 6-11 = 0$; similarly, $\text{wt}_f(3) = 0$ and  $\text{wt}_f(4)=0$. 
    With that we have verified that $2\ar{C_6}$ is a $3-$node ornament, which is not an isolated case. It shows that $2\ar{C_L}$ is a $3-$node ornament given any even  $L\geq 6$.  
\begin{theorem} \label{thm:3nodes2cycles}
     Let $L=2k$ 
     for any integer $k\ge 3$. Then there exists a labeling $h$ on $2\ar{C_L}$ such that $2\ar{C_L}$ is a $3-$nodes ornament. Thus, if   $\ar{G_0}$ is a DDMOG with $n\ge 4$ 
     vertices, then $\ar{G_0} \oplus 2\ar{C_L}$ is a DDMOG.  
\end{theorem}

\begin{figure}
\centering

\begin{tikzpicture}[scale=0.67, transform shape]
  \tikzset{
    labeled_vertex/.style = {circle, draw=black, fill=white, inner sep=1.5pt, text=black, minimum size=24pt, very thick},
    arc/.style = {->, >=stealth, very thick},
    dots_node/.style = {inner sep=2pt, minimum size=24pt, text=black},
    rlab/.style = {text=red!85!black, font=\small\bfseries}
  }

  \node [labeled_vertex, label={[rlab]left:$+1$}] (1) at (-10, 1) {$v_1$};
  \node [labeled_vertex, label={[rlab]above left:$-3$}] (2) at (-8, 2) {$v_2$};
  \node [labeled_vertex, label={[rlab]below left:$+4$}] (3) at (-8, 0) {$v_3$};
  \node [labeled_vertex, label={[rlab]above:$-4$}] (4) at (-5.5, 2) {$v_4$};
  \node [labeled_vertex, label={[rlab]below:$+4$}] (5) at (-5.5, 0) {$v_5$};
  
  \node [dots_node] (dd11) at (-3, 2) {$\dots$};
  \node [dots_node] (dd14) at (-3, 0) {$\dots$};
  
  \node [labeled_vertex, label={[rlab]above:$-4$}] (i-2) at (-0.5, 2) {$v_{i-2}$};
  \node [labeled_vertex, label={[rlab]below:$+4$}] (i-1) at (-0.5, 0) {$v_{i-1}$};
  \node [labeled_vertex, label={[rlab]above:$-4$}] (i) at (2, 2) {$v_i$};
  \node [labeled_vertex, label={[rlab]below:$+4$}] (i+1) at (2, 0) {$v_{i+1}$};
  \node [labeled_vertex, label={[rlab]above:$-4$}] (i+2) at (4.5, 2) {$v_{i+2}$};
  \node [labeled_vertex, label={[rlab]below:$+4$}] (i+3) at (4.5, 0) {$v_{i+3}$};
  
  \node [dots_node] (dd12) at (7, 2) {$\dots$};
  \node [dots_node] (dd13) at (7, 0) {$\dots$};
  
  \node [labeled_vertex, label={[rlab]above right:$-4$}] (L-2) at (9.5, 2) {$v_{L-2}$};
  \node [labeled_vertex, label={[rlab]below right:$+3$}] (L-1) at (9.5, 0) {$v_{L-1}$};
  \node [labeled_vertex, label={[rlab]right:$-1$}] (L) at (11.5, 1) {$v_L$};

  \draw[arc] (1) to (2);
  \draw[arc] (2) to (4);
  \draw[arc] (4) to (dd11);
  \draw[arc] (dd11) to (i-2);
  \draw[arc] (i-2) to (i);
  \draw[arc] (i) to (i+2);
  \draw[arc] (i+2) to (dd12);
  \draw[arc] (dd12) to (L-2);
  \draw[arc] (L-2) to (L);
  \draw[arc] (L) to (L-1);
  \draw[arc] (L-1) to (dd13);
  \draw[arc] (dd13) to (i+3);
  \draw[arc] (i+3) to (i+1);
  \draw[arc] (i+1) to (i-1);
  \draw[arc] (i-1) to (dd14);
  \draw[arc] (dd14) to (5);
  \draw[arc] (5) to (3);
  \draw[arc] (3) to (1);

  \node [labeled_vertex, label={[rlab]left:$-1$}] (1') at (-10, -5) {$v_1'$};
  \node [labeled_vertex, label={[rlab]above left:$+3$}] (2') at (-8, -4) {$v_2'$};
  \node [labeled_vertex, label={[rlab]below left:$-4$}] (3') at (-8, -6) {$v_3'$};
  \node [labeled_vertex, label={[rlab]above:$+4$}] (4') at (-5.5, -4) {$v_4'$};
  \node [labeled_vertex, label={[rlab]below:$-4$}] (5') at (-5.5, -6) {$v_5'$};
  
  \node [dots_node] (dd21) at (-3, -4) {$\dots$};
  \node [dots_node] (dd24) at (-3, -6) {$\dots$};
  
  \node [labeled_vertex, label={[rlab]above:$+4$}] (i-2') at (-0.5, -4) {$v_{i-2}'$};
  \node [labeled_vertex, label={[rlab]below:$-4$}] (i-1') at (-0.5, -6) {$v_{i-1}'$};
  \node [labeled_vertex, label={[rlab]above:$+4$}] (i') at (2, -4) {$v_i'$};
  \node [labeled_vertex, label={[rlab]below:$-4$}] (i+1') at (2, -6) {$v_{i+1}'$};
  \node [labeled_vertex, label={[rlab]above:$+4$}] (i+2') at (4.5, -4) {$v_{i+2}'$};
  \node [labeled_vertex, label={[rlab]below:$-4$}] (i+3') at (4.5, -6) {$v_{i+3}'$};
  
  \node [dots_node] (dd22) at (7, -4) {$\dots$};
  \node [dots_node] (dd23) at (7, -6) {$\dots$};
  
  \node [labeled_vertex, label={[rlab]above right:$+4$}] (L-2') at (9.5, -4) {$v_{L-2}'$};
  \node [labeled_vertex, label={[rlab]below right:$-3$}] (L-1') at (9.5, -6) {$v_{L-1}'$};
  \node [labeled_vertex, label={[rlab]right:$+1$}] (L') at (11.5, -5) {$v_L'$};

  \draw[arc] (2') to (1');
  \draw[arc] (1') to (3');
  \draw[arc] (3') to (5');
  \draw[arc] (5') to (dd24);
  \draw[arc] (dd24) to (i-1');
  \draw[arc] (i-1') to (i+1');
  \draw[arc] (i+1') to (i+3');
  \draw[arc] (i+3') to (dd23);
  \draw[arc] (dd23) to (L-1');
  \draw[arc] (L-1') to (L');
  \draw[arc] (L') to (L-2');
  \draw[arc] (L-2') to (dd22);
  \draw[arc] (dd22) to (i+2');
  \draw[arc] (i+2') to (i');
  \draw[arc] (i') to (i-2');
  \draw[arc] (i-2') to (dd21);
  \draw[arc] (dd21) to (4');
  \draw[arc] (4') to (2');

\end{tikzpicture}

\caption{Generalized $3-$node ornament, $2\ar{C_L}$, with labeling $h(v_k)=n+k$ and $h(v'_k)=n+L+k$ for $1\le k\le L$, given any $n\ge4$. The values of wt$_h$ are marked in red next to each vertex. }
\label{fig:2GenNgons3nodes}
\end{figure}
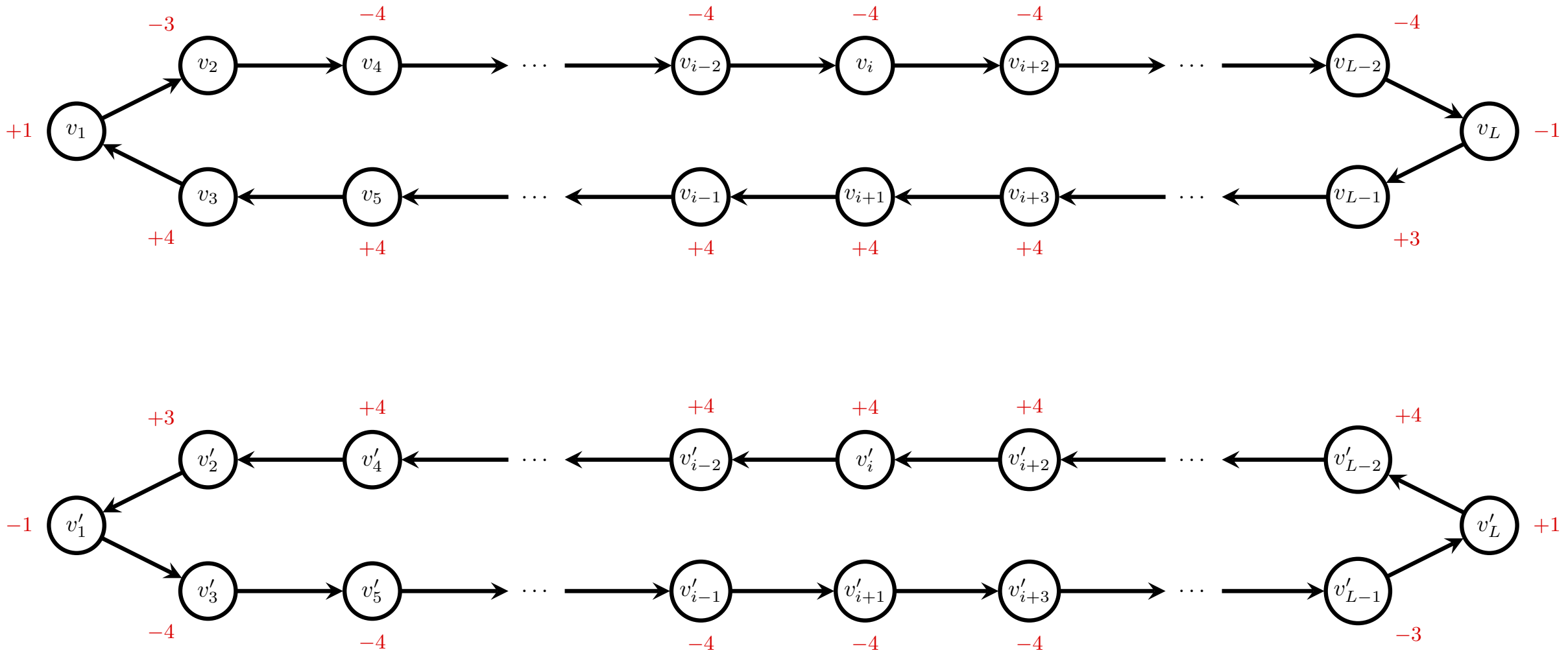
\begin{proof}
    Let $L=2k$, $k\geq 3$. Let $\ar{G_0}$  be an oriented graph with $n\ge 4$ vertices $v_1, v_2, \hdots, v_n$, and a DDM labeling bijection $f_0:V(G_0)\rightarrow \{1,...,n\}$, defined as $f_0(v_p)=p$ for any $p \in \{1,2,\hdots, n\}$.  For $i \in \{1, 2\}$, let 
    $H_i$ be an $L-$cycle with  vertex sets  $V(H_{1})= \{v_q \}_{q=1}^L$, $V(H_{2})= \{v'_q \}_{q=1}^L$.
       We now define an orientation on $H_1$ and $H_2$ by defining the  edges sets  as follows: 
      \begin{eqnarray}
          E(\ar{H_1})&=\{v_1v_3,v_3v_5,v_5v_7,...,v_{L-1}v_L,v_Lv_{L-2},v_{L-2}v_{L-4},...,v_4v_2,v_2v_1\}, \\
          E(\ar{H_2}) &=\{v'_3v'_1,v'_5v'_3,v'_7v'_5,...,v'_{L}v'_{L-1},v'_{L-2}v'_{L},v'_{L-4}v'_{L-2},...,v'_2v'_4,v'_1v'_2\}.
      \end{eqnarray}

\noindent 
Let $H=H_1\cup H_2$, and let  $h: V(H) \rightarrow \{n+1, ..., n+2L\}$ be the bijection defined for any $q \in \{1,2,\hdots, L\}$ by $h(v_q)=n+q$ when $v_q \in V(H_1)$,  and $h(v'_q)=n+L+q $ when $ v_q' \in V(H_2)$.
Observe   that the weight $\text{wt}_h$ at any vertex in $ {H}$ has values $\pm1$, $\pm3$, or $\pm4$ (see Figure \ref{fig:2GenNgons3nodes}). Notice that, if $\ar{H}= \ar{H_1} \cup \ar{H_2}$ is an ornament, it would be a $3-$nodes ornament, and 
$\displaystyle \max_{v\in V(H)}\{\vert \text{wt}_h(v)\vert \} = \max\{1, 3, 4\} =4 <5 \leq n$. 
It is easy to see that $\displaystyle \sum_{v \in V_h ^1(H)}h(v)=h(v_1)+h(v_{L}')=(n+1)+(n+L+1)=2n+L+2$
 equals $\displaystyle \sum_{v \in V_h ^{-1}(H)}h(v)=h(v_1')+h(v_{L})=(n+L+1)+(n+L)=2n+L+2$. 
Similarly,    $\displaystyle \sum_{v \in V_h ^3(H)}h(v) = 2n+2L+1=\sum_{v \in V_h ^{-3}(H)}h(v)$. Finally, it can be verified that 

\noindent $\displaystyle \sum_{v \in V_h ^4(H)}h(v) = \sum_{j=2}^{\frac{L-4}{2}} h(v_{2j-1}) + \sum_{j=2}^{\frac{L-2}{2}} h(v_{2j}')= nk(k-2)(k-4) = \sum_{j=2}^{\frac{L-2}{2}} h(v_{2j}) +   \sum_{j=2}^{\frac{L-4}{2}} h(v_{2j-1}') =\displaystyle \sum_{v \in V_h ^{-4}(H)}h(v)$.

 We define a bijective labeling $f$ for the new graph $\ar{G}=\ar{G_0}\oplus_h \cup_{i=1}^2 \ar{H_i}$  as in \eqref{deffonG0H}, with  $\ar{H}=\bigcup_{i=1}^2 \ar{H_i}$.   
  The elements of the set of additional edges $E'$ is defined as follows: 
For any vertex $v \in V(H)$ with weight $wt_h(v)=p \in \{1,3,4\}$, we introduce an edge $vu_p$ going from $v$ to the vertex $u_p\in V(G_0)$ with the label $p$. Likewise, if the weight of $v$ is $wt_h(v)=-p$, we introduce the edge $u_pv$, where $f_0(u_p)=p$.
  With that, by Theorem~\ref{thmWgt1} we have that $\ar{G}$ is a DDMOG; thus $\cup_{i=1}^2 \ar{H_i}$ is an ornament.
\end{proof}

\begin{remark}
    The labeling  used in the proof of Theorem~\ref{thm:3nodes2cycles}  (as illustrated in  Figure \ref{fig:2GenNgons3nodes})  can also be employed to create an ornament when $L=2k+1$ for any integer $k\ge 2$, but in  this case the ornament must consist of four cycles. We leave the verification of this statement  as an exercise for the reader.
\end{remark}

\begin{example}\label{ex:2c43nodes}  $2C_L$ can be an ornament for $L<6$ even though Theorem~\ref{thm:3nodes2cycles} only offers a solution for $L\geq 6$. Specifically,  there exists a labeling on $2\ar{C_4}$ such that it is a $2-$node ornament.
To verify this, consider first an existing DDMOG $\ar{G_0}$ with $n$ vertices.
Here's a brief argument  why $\ar{G}=\ar{G_0}\oplus 2\ar{C_4}$ is a DDMOG: Let  $f$ denote the DDM labeling on $\ar{G}$; we  analyze the weight induced by   $f$.  The weights at the nodes in  $ \ar{G_0}$, labeled $1$ and $3$, are as follows  (Subfigure \ref{fig:thm1L4}):

\noindent wt$_f(1)=(n+1)+(n+8)-(n+4)-(n+5)=0$, 
wt$_f(3)= (n+3) + (n+6)-(n+2)-(nj+5)=0$.

Thus, only $2\ar{C_4}$ are required to create a $2-$nodes ornament. 
\end{example}

  \subsection{Minimal value for $L$ when $4\ar{C_L}$ is a 3-node ornament  }\label{subs:minLornament} 
  Theorem~\ref{theorem:4CL-ornaments} guarantees the construction of   $3-$nodes ornaments   for $L\geq 6$, but Example~\ref{ex:2c43nodes}  shows that the minimal value of $L$ is at most $4$. While the labeling pattern used in the proof of Theorem~\ref{theorem:4CL-ornaments} is repetitive within each used cycle,  abandoning that constraing let us to find a pattern  for $4\ar{C_3}$ as a $3-$nodes ornament. Thus,  we can state  the minimal value for $L$ when  $4\ar{C_L}$ is a $3-$nodes ornament: 
  \begin{corollary}
      $\min\{L \in \mathbb{Z}^+\; | \; \text{$4\ar{C_L}$ is a $3-$nodes ornament }\} = 3. $
  \end{corollary}
  \begin{proof}
      Clearly,  $\min\{L \in \mathbb{Z}^+\; | \; \text{$4\ar{C_L}$ is a $3-$nodes ornament }\} \geq 3. $ Consider the labeling $h$ as defined in Figure~\ref{3nodesC3}. Working with any oriented graph $\ar{G_0} $ with $n\geq 4$ vertices and a DDM labeling $f_0$,   we have $V_h^1 = \{v_6, v_7\}$ and $V_h^{-1} = \{v_1, v_{12}\}$; thus, $\sum_{v \in V_h^1} h(v)= h(v_6) + h(v_7)= n+6 + n +7 = n+1 +n+12 =h(v_1) + h(v_{12})= \sum_{v \in V_h^{-1}} h(v)$. Similarly,  $\sum_{v \in V_h^a}  = \sum_{v \in V_h^{-a}}$ for $a \in \{3, 4\}$. By  Theorem~\ref{thmWgt1}, $G_0\oplus_{\text{wt}_h} (4\ar{C_L})$ is a DDMOG. Thus, $4\ar{C_L}$ is an ornament, with attachment nodes labeled $1$, $3$, $4$. 
  \end{proof}

\begin{figure}[h]
\centering

\begin{tikzpicture}[scale=0.76]
    \tikzset{
        labeled_vertex/.style = {circle, draw=black, fill=white, inner sep = 1.5pt, text=black, minimum size=20pt, very thick},
        arc/.style={->,very thick}
    };

    \node [labeled_vertex, label={[text=red]left:-1}] (1) at (3, -1) {\(v_1\)};
    \node [labeled_vertex, label={[text=red]above:+4}] (4) at (5, -2) {\(v_4\)};
    \node [labeled_vertex, label={[text=red]left:-3}] (5) at (3, -3) {\(v_5\)};

    \node [labeled_vertex, label={[text=red]left:+1}] (6) at (7, -1) {\(v_6\)};
    \node [labeled_vertex, label={[text=red]above:-4}] (3) at (9, -2) {\(v_3\)};
    \node [labeled_vertex, label={[text=red]left:+3}] (2) at (7, -3) {\(v_2\)};

    \node [labeled_vertex, label={[text=red]left:+1}] (7) at (11, -1) {\(v_7\)};
    \node [labeled_vertex, label={[text=red]left:+3}] (11) at (11, -3) {\(v_{11}\)};
    \node [labeled_vertex, label={[text=red]above:-4}] (10) at (13, -2) {\(v_{10}\)};

    \node [labeled_vertex, label={[text=red]left:-1}] (12) at (15, -1) {\(v_{12}\)};
    \node [labeled_vertex, label={[text=red]left:-3}] (8) at (15, -3) {\(v_8\)};
    \node [labeled_vertex, label={[text=red]above:+4}] (9) at (17, -2) {\(v_9\)};

    \draw[arc] (1) to (5);
    \draw[arc] (5) to (4);
    \draw[arc] (4) to (1);

    \draw[arc] (6) to (2);
    \draw[arc] (2) to (3);
    \draw[arc] (3) to (6);

    \draw[arc] (7) to (10);
    \draw[arc] (10) to (11);
    \draw[arc] (11) to (7);

    \draw[arc] (9) to (8);
    \draw[arc] (8) to (12);
    \draw[arc] (12) to (9);
\end{tikzpicture}
 
\caption{$3-$nodes ornament $4\ar{C_3}$, with labeling defined as follows: $h(v_i)=n+i$, $1\leq i \leq 12$; here $n\geq 4$. The values of $\text{wt}_h (v)$ are marked in red next to each vertex $v$.  }
\label{3nodesC3}
\end{figure}
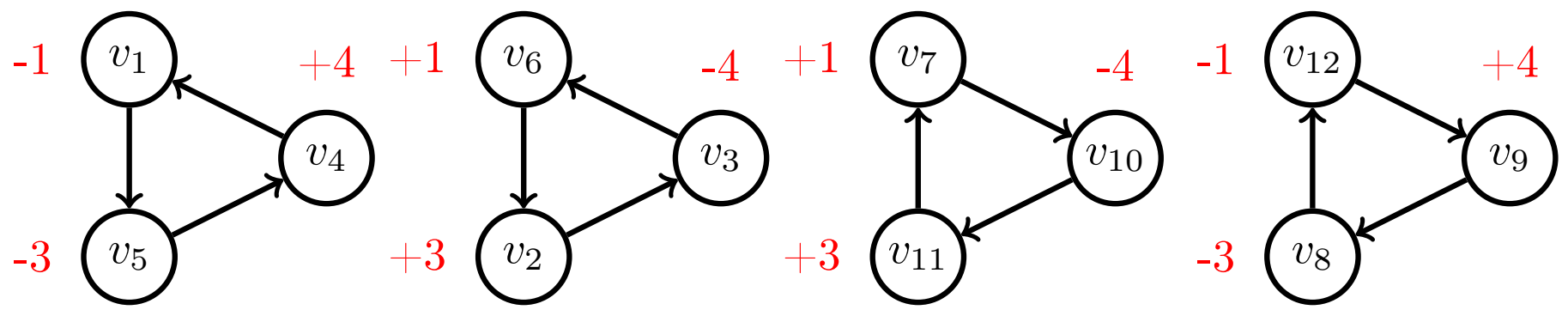

\subsection{{A} $2-$nodes ornament, formed by one $L-$cycle }  \label{subs:onecycle}
For special values of $L$, it shows that $\ar{C_L}$ is an ornament (with a carefully selected labeling function).
\begin{theorem}\label{thm:2nodes1cycle}
    Let $L=4k$ for some integer $k\geq 2$. 
    There exists a labeling $h$ on the vertices of $H=C_{L}$ and there exists an orientation on $H$ such that  
    $\ar{H}$ is a $2-$nodes ornament   for all $n\geq 4$, with nodes of attachment labeled by $2$ and $4$.
\end{theorem}

\begin{proof} Let $\ar{G_0}$ be a DDMOG with $n\geq 4$ vertices, paired with a DDM labeling $f_0$. Let $u_2$,  $u_4$ $\in V(G_0)$   such that $f_0(u_2)=2$, $f_0(u_4)=4$. Let $\ar{H}=\ar{C_{4k}}$ for some $k\geq 2$, with the orientation and labeling  as 
in Figure~\ref{2nodeAttach1}.

    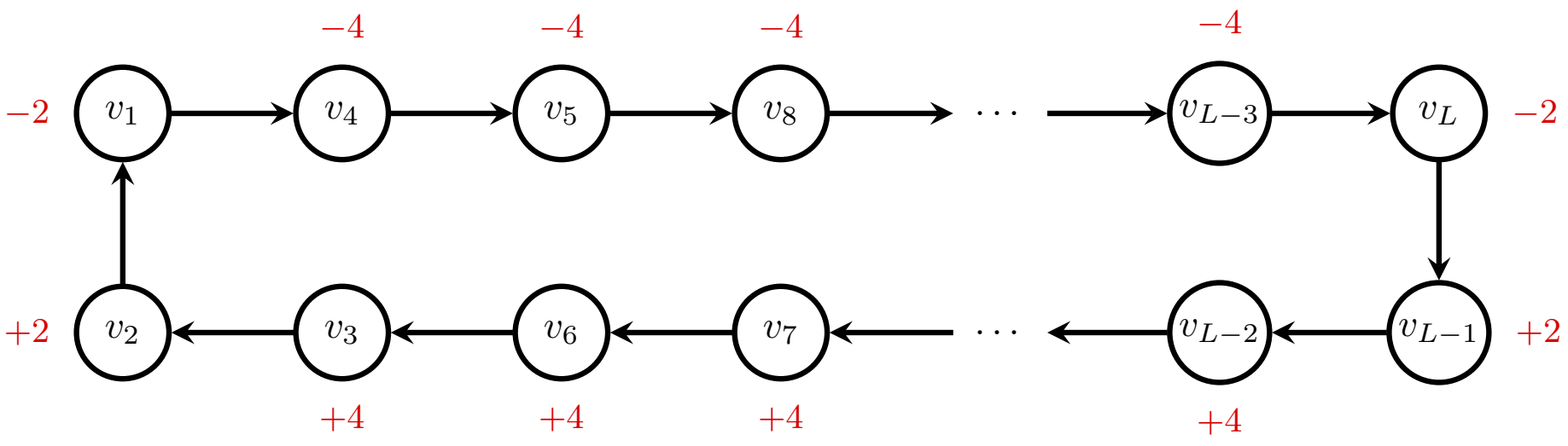
\begin{figure}[h]
\centering

\begin{tikzpicture}[scale=0.85, transform shape]
  \tikzset{
    labeled_vertex/.style = {circle, draw=black, fill=white, inner sep=1.5pt, text=black, minimum size=24pt, very thick},
    arc/.style = {->, >=stealth, very thick},
    dots_node/.style = {inner sep=2pt, minimum size=24pt, text=black},
    rlab/.style = {text=red!85!black, font=\small\bfseries}
  }

  \node [labeled_vertex, label={[rlab]left:$-2$}] (1) at (-7, 2) {$v_1$};
  \node [labeled_vertex, label={[rlab]above:$-4$}] (4) at (-5, 2) {$v_4$};
  \node [labeled_vertex, label={[rlab]above:$-4$}] (5) at (-3, 2) {$v_5$};
  \node [labeled_vertex, label={[rlab]above:$-4$}] (8) at (-1, 2) {$v_8$};
  
  \node [dots_node] (dd) at (1, 2) {$\dots$};
  
  \node [labeled_vertex, label={[rlab]above:$-4$}] (4k-3) at (3, 2) {$v_{L-3}$};
  \node [labeled_vertex, label={[rlab]right:$-2$}] (4k) at (5, 2) {$v_L$};

  \node [labeled_vertex, label={[rlab]left:$+2$}] (2) at (-7, 0) {$v_2$};
  \node [labeled_vertex, label={[rlab]below:$+4$}] (3) at (-5, 0) {$v_3$};
  \node [labeled_vertex, label={[rlab]below:$+4$}] (6) at (-3, 0) {$v_6$};
  \node [labeled_vertex, label={[rlab]below:$+4$}] (7) at (-1, 0) {$v_7$};
  
  \node [dots_node] (dd2) at (1, 0) {$\dots$};
  
  \node [labeled_vertex, label={[rlab]below:$+4$}] (4k-2) at (3, 0) {$v_{L-2}$};
  \node [labeled_vertex, label={[rlab]right:$+2$}] (4k-1) at (5, 0) {$v_{L-1}$};

  \draw[arc] (2) to (1);
  \draw[arc] (1) to (4);
  \draw[arc] (4) to (5);
  \draw[arc] (5) to (8);
  \draw[arc] (8) to (dd);
  \draw[arc] (dd) to (4k-3);
  \draw[arc] (4k-3) to (4k);
  \draw[arc] (4k) to (4k-1);
  \draw[arc] (4k-1) to (4k-2);
  \draw[arc] (4k-2) to (dd2);
  \draw[arc] (dd2) to (7);
  \draw[arc] (7) to (6);
  \draw[arc] (6) to (3);
  \draw[arc] (3) to (2);

\end{tikzpicture}

\caption{$2-$node ornament, formed by one $(4k)-$cycle,  $k\geq 1$. The  labeling  is  defined as follows: $h(v_j)=n+j$, for $1\leq j \leq 4k$; here $n\geq4$. The values of $\text{wt}_h (v)$ are marked in red next to each vertex $v$. }
\label{2nodeAttach1}
\end{figure}
Observe that the weight values $\text{wt}_h$ at every vertex $v\in V(H)$ are either $\pm 2$ or $\pm 4$; thus if $\ar{H}$ is an ornament, it would be a $2-$nodes ornament, and 
\begin{equation} \label{t1test1}
\max_{v\in V(H)}\{\vert \text{wt}_h(v)\vert \} = \max\{2, 4\} =4 <5 \leq n. \end{equation}
Note that $\displaystyle \sum_{v \in V_h ^2(H)}h(v)=h(v_2)+h(v_{L-1})=(n+2)+(n+L-1)=2n+L+1= (n+1)+ (n+L)=h(v_1)+h(v_{L})= \sum_{v \in V_h ^{-2}(H)}h(v)$. In a similar fashion, we can verify that 

\noindent $\displaystyle \sum_{v \in V_h ^4(H)}h(v)=h(v_3) +h(v_6) + h(v_7)+ \hdots +h(v_{L-2})= (\frac{L}{4}-1)(2n+L+1) =\displaystyle \sum_{v \in V_h ^{-4}(H)}h(v)$. Thus, 
\begin{equation}  \label{t1test2}\sum_{v \in V_h ^{-i}(H)}h(v)=\sum_{v \in V_h ^4(H)}h(v) \; \text{ for $i \in \{2, 4\}$.}\end{equation}
With the labeling $f$ defined as in \eqref{deffonG0H}, and by \eqref{t1test1}, \eqref{t1test2} and Theorem~\ref{thmWgt1}, it follows that $\ar{G}=\ar{G_0} \oplus_h \ar{C_{4k}}$  is a DDMOG, therefore, $\ar{H}$ is an ornament.

\end{proof}
The new DDMOG  $\ar{G}=\ar{G_0} \oplus_h \ar{C_{4k}}$  as per the construction in Theorem~\ref{thm:2nodes1cycle} has an edges set $E(G)=E(G_0)\cup E(H) \cup E'$, where $E'$ contains new oriented edges between $u_2$, $u_4$ $\in V(G_0)$ on one side, and the vertices in $ V(H)$ on the other. Specifically, for every  $v\in V(H)$, and $p \in \{2, 4\}$ we have 
\begin{center}
i.   if $\displaystyle \text{w}_h (v) =-p$, then $\displaystyle u_pv \in E'$; \hspace{13mm} ii.   if $\displaystyle \text{w}_h (v) =p$, then $\displaystyle vu_p \in E'$.
\end{center}
The DDM labeling $f$ on $V(\ar{G})$ is defined in the usual way, as in equation \eqref{deffonG0H}.  

\begin{remark}\label{rem:specialK=1}
    In the special case when $k=1$, the labeling used in the proof of Theorem~\ref{thm:2nodes1cycle} (as seen in Figure~\ref{2nodeAttach1}) results with a  $1-$node ornament   formed by a single cycle $\ar{C_4}$. 
\end{remark}

\subsection{Variations in labeling $p\ar{C_4}$, $p\geq 1$}
\label{section:minimal-L}
\label{squareO}

A couple of  different ways to label an ornament of type $p\ar{C_4}$, $p\geq 1$, were already discussed in \cite{FM, ALR}. This section expands on those labelings, and offers other ways  to create ornaments that involve less than four cycles   or allow for an increased number of node attachments. We also offer a result that  highlights the special property of $4-$cycles. Namely, if $p\ar{C_L}$ is  a $1-$node ornament, then $L$ must be equal to $4$: 

\begin{proposition}
\label{prop:c4is1nodeonly}    
 Let $p\geq 1$, and let   $L\geq 3$ but  $L\neq 4$. Then  $p\ar{C_L}$ cannot be a $1-$node ornament.\end{proposition}
 \begin{proof} We will use the method of contradiction. Let $k >0$ be a  fixed integer. Suppose   $p\ar{C_L}$ is a $1-$node ornament for some fixed $p\geq 1$, with associated labeling $h$ which induces a weight $\text{wt}_h$ with values $k$ or $-k$ at every vertex. Under these assumptions, we will check what the constraints are on the values of $L$. 
 \begin{itemize}
 \item[i.] Let  $L \geq 3$ be an odd integer.  
Thus within one cycle $\ar{C_L}$ with vertices $\{v_i\}_{i=1}^L$, we have:
\begin{center} $h(v_{i-1}) - h(v_{i+1}) = \pm k$, $2\leq i \leq L-1$,   $h(v_{L-1})-h(v_1)= \pm k$, and
  $h(v_L)-h(v_2)=\pm k$.\end{center}
\noindent By taking the sum of the first $L-1$ equations we obtain that $h(v_L)-h(v_2) = tk$ for some $t\in \mathbb{Z}$ which should equal $\pm k$ by the last equation. Since $L$ is assumed odd, we sum an even number of times values of $k$ or $-k$, which could equal $\pm 2mk$ for some $m \in \{0, 1, \hdots, (L-1)/2\}$, therefore will never equal $\pm k$. Thus,  $L$ cannot be odd.

\item[ii.] Let $L=4m+2$, $m\geq 1$. Observe that within one cycle we have
\begin{center}
    $h(v_{2i-1})-h(v_{2i+1})= \pm k$, $1\leq i\leq 2m$, $h(v_{4m+1}) - h(v_1)= \pm k$. 
\end{center}
\noindent 
Then on one hand we have $\displaystyle \sum_{i=1}^{2m} h(v_{2i-1})-h(v_{2i+1}) =0$, but that sum must equal a sum of 
an odd number of addends, each with values of either $k$ and $-k$. This sum  cannot be $0$ no matter how many positive versus negative values we have; thus, $L$ cannot take values of type $4m+2$, $m\geq 1$.

\item[iii.] Let $L=4m$, $m\geq 2$.   Within one oriented cycle, for the odd indices  we have: 
\begin{equation}\label{odddiffk}
    h(v_{2i-1}) - h(v_{2i+1}) = \pm k, \; 1\leq i \leq 2m-1, \; \text{and} \; h(v_{4m-1}) - h(v_1) = \pm k.
\end{equation}
If we add up all of the $2m$ equations in \eqref{odddiffk}, we get zero on the left hand side, as all $h(v_j)$ cancel out; given that on the right hand side we have addends with values either  $- k$ or $k$, it must be that exactly $m $ equations have $k$ as an output, and $m$ equations have $- k$ as an output. Now if we only add the first two equations, we get 
$h(v_1) -h(v_5) = \pm 2tm$, where $t$ is either $0$ or $1$. If $t=0$ then $h(v_1)=h(v_5)$ which cannot happen, so it must be that $t=1$, that is, the first two equations have the same output (either $k$ or $-k$). Similarly, by adding two consecutive equations, we conclude that all of the equations must have the same  output, which contradicts the previous conclusion that exactly $m $ equations have $k$ as an output, and $m$ equations have $- k$ as an output. We must conclude that our initial assumption is wrong.
\end{itemize}
Therefore, $4\ar{C_{L}}$ cannot be a $1-$node ornament for any $L\geq 3$ but $L\neq 4$. 
\end{proof}

\subsubsection{Alternative labeling of a $1-$node $4-$cycle ornament}~\label{4cycle2notsPart1}
We present here a new  labeling on $p\ar{C_4} $   to create a $1-$node ornament that differs from the one introduced in \cite{AAetc}. This new labeling  involves uniform jumps between the labels of vertices in each cycle, and results in an attachment to a single node $u$ with label $j$, where $j\geq 1$ is the jump in the labeling.  Note that the aforementioned  labeling  in \cite{AAetc} only allowed for an attachment node labeled with $2j$, $j\geq 1$.

We organize the labeling function $h$ and the weights $\text{wt}_h$ on a $4-$cycle in Table~\ref{table:j-4-cycles-ornament}.
  \begin{table}
  \centering
\begin{tabular}{ c | c | c  }
Vertex, $v$ & $h(v)$ & wt$_h$($v$) \\
\hline
$v_{i1}$ & $n+i$ & $(n+3j+i)-(n+2j+i)=j $ \\ 
$v_{i2}$ & $n+j+i$ & $(n+2j+i)-(n+3j+i)=-j  $ \\  
 $v_{i3}$ & $n+2j+i$ &  $(n+i)-(n+j+i)=-j$ \\
 $v_{i4}$ & $n+3j+i$ & $ (n+j+i)-(n+i)=j$ 
\end{tabular}
\caption{The labeling $h$ of \(jC_4\) as a \(1\)-node ornament, attaching to a node with label $j$, $j\geq 1$.  }
\label{table:j-4-cycles-ornament}
\end{table}
To create a new DDMOG with this labeling we  need a DDMOG $\ar{G_0}$ with $n$ vertices, where $n\geq j$. We  also need $p=j$ many $4-$cycles $\ar{H}_i$, $1\leq i \leq j$. 
Note that before introducing new edges from or to $u_j \in V(\ar{G_0})$ such that $f_0(u_j)=j$,   the weights at the vertices of $\cup_{i=1}^j\ar{H}_i$  are $\pm j$. We introduce the new edges   as follows: $ (v_{i1}, u_j)$, $(u_j, v_{i2})$, $(u_j, v_{i3})$, $(v_{i4}, u_j)$. With these edges, the weight at all vertices in $V(\cup_{i=1}^j\ar{H}_i)$ is 0. The weight at the node $u_j$ is also 0 since the only change to the weight was $+2j-2j=0$ per cycle. Lastly, in order for these cycles to create an ornament the orientation on all $j$ $4-$cycles must be the same. 
Therefore, we have the following result: 

\begin{corollary}
    \label{corollary:jC_4-is-1-node-ornament}
    Given any $j\geq 1$, the oriented graph $j\ar{C_4}$ is a $1-$node ornament, with node of attachment labeled by $j$. 
\end{corollary}

\begin{remark}\label{onenodeCycleM} The wheel  $\ar{W_4}$ (Figure~\ref{fig:a}) can be seen as a weighted sum of the one-vertex graph $\ar{K_1}$ with the $4-$cycle $\ar{C_4}$. The windmill of wheels $\ar{W_4}^{(k)}$,  introduced in \cite{AAetc}, can also be created by attaching $k$ copies of $\ar{C_4}$ to $\ar{K_1}$. We highlight that  this new labeling as presented in  Table~\ref{table:j-4-cycles-ornament},  generates a DDM labeling for   $\ar{W_4}^{(k)}$, without the need to change the labeling for every additional cycle being added to the  ``evergrowing'' windmill. The construction in \cite{AAetc} unfortunately requires complete relabeling with every new cycle added.  \end{remark}
%
%
\subsubsection{Minimal value for $L$ in a $s-$nodes ornament when $s \in \{4, 6\}$ }\label{subsection6C4}
  In  \cite{ALR} the authors  have shown that $4\ar{C_L}$ is a  $6-$nodes ornament when $L\geq 8$; that is, the minimal known value of $L$ so that $4\ar{C_L}$ is a $6-$nodes ornament was $L=8$. The proof of that result was based on a repetitive labeling pattern within  every copy of the $L-$cycle.   If we increase the number of cycles used and relax the request for repetitive labeling as well as for the number of cycles used, it is possible to create a $6-$nodes ornament from $4-$cycles.  
  
  \begin{proposition}
      There exists a labeling $h$ on $\ar{H}=6\ar{C_4}$ such that $\ar{H}$ is a $6-$nodes ornament.
  \end{proposition}
  \begin{proof}
We will show that the labeling $h$ on  $\ar{H}=6\ar{C_4}$ as presented in  Figure~\ref{fig:6nodes4vertex} allows for  $\ar{H}$ to  be attached to any DDMOG $\ar{G_0}$ with at least $11$ vertices in a way that the resulting graph is DDMO. Suppose  the DDM labeling on $\ar{G_0}$ is $f_0$, with   $f_0(u_p)=p$ for any $u_p \in V(G_0)$.  Then the nodes are those $u_p \in V(\ar{G_0})$, where $p\in N(\ar{H})= \{1,3,5,7,9,11\}$.  
 The new DDM labeling  introduced on the vertices of $\ar{G}:=\ar{G_0}\oplus_h  \ar{H}$ will be defined as in \eqref{deffonG0H}. 
 At any $p\in N(\ar{H})$, if the weight induced by the labeling $h$ as presented in the figure is $\text{wt}_h(v)=p$ at a  vertex $v \in V( \ar{H})$, then a new edge is introduced, $vu_p$; if $\text{wt}_h(v)=-p<0$ at  $v \in V(\ar{H})$, then a new edge is introduced, $u_pv$. It is easy to verify that these additional edges ensure that $\text{wt}_f$ is a zero function at every vertex in $V( \ar{G})$. Specifically, at the nodes we have  $\text{wt}_f (u_p) =0$ for all $p \in N(\ar{H})$; for instance:  
 $\text{wt}_f (u_1) = (n+8)-(n+5) + (n+17)-(n+20)=0$. Thus, $\ar{G} $ is a DDMOG, and $\ar{H}=6\ar{C_4}$ is a $6-$nodes ornament.
  \end{proof} 



\begin{example}\label{ex:2nodes2C4}
It is known \cite{ALR} that  $4\ar{C_4}$ is  a $2-$nodes ornament; one labeling is seen in Subfigure~\ref{fig:thm1L4}, where the labels of the attachment nodes need to be $1$ and $3$. By rearranging the labels, it is   easy to show that the labels of the attachment nodes can be $2$ and $6$ (Subfigure~\ref{fig:thm1L4b}). 
\end{example}

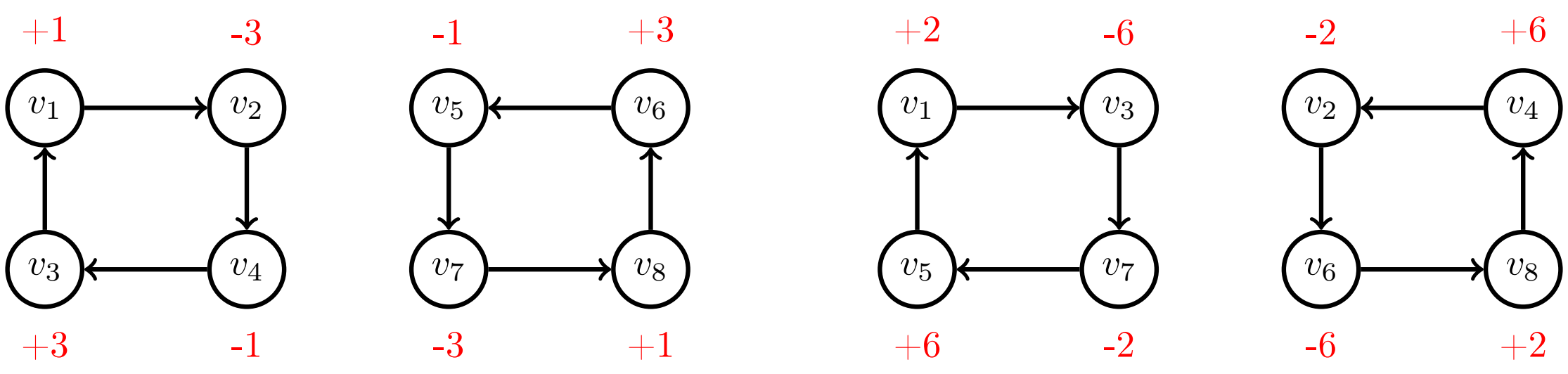
\begin{figure}[htbp]
\centering

\begin{subfigure}[b]{0.49\textwidth}
\centering
\begin{tikzpicture}[scale=0.76]
    \tikzset{
        labeled_vertex/.style = {circle, draw=black, fill=white, inner sep = 1.5pt, text=black, minimum size=20pt, very thick},
        arc/.style={->,very thick}
    };

    \node [labeled_vertex, label={[text=red]above:+1}] (4) at (-2, 2) {\(v_1\)};
    \node [labeled_vertex, label={[text=red]above:-3}] (5) at (0.5, 2) {\(v_2\)};
    \node [labeled_vertex, label={[text=red]below:-1}] (7) at (0.5, 0) {\(v_4\)};
    \node [labeled_vertex, label={[text=red]below:+3}] (6) at (-2, 0) {\(v_3\)};

    \node [labeled_vertex, label={[text=red]above:-1}] (8) at (3, 2) {\(v_5\)};
    \node [labeled_vertex, label={[text=red]above:+3}] (9) at (5.5, 2) {\(v_6\)};
    \node [labeled_vertex, label={[text=red]below:+1}] (11) at (5.5, 0) {\(v_8\)};
    \node [labeled_vertex, label={[text=red]below:-3}] (10) at (3, 0) {\(v_7\)};

    \draw[arc] (4) to (5);
    \draw[arc] (5) to (7);
    \draw[arc] (7) to (6);
    \draw[arc] (6) to (4);

    \draw[arc] (10) to (11);
    \draw[arc] (11) to (9);
    \draw[arc] (9) to (8);
    \draw[arc] (8) to (10);
\end{tikzpicture}
 
\caption{ $2\ar{C_4}$, with  nodes labeled $1$ and $3$.}
\label{fig:thm1L4}
\end{subfigure}
%
%
\begin{subfigure}[b]{0.49\textwidth}
\centering
 
\begin{tikzpicture}[scale=0.76]
    \tikzset{
        labeled_vertex/.style = {circle, draw=black, fill=white, inner sep = 1.5pt, text=black, minimum size=20pt, very thick},
        arc/.style={->,very thick}
    };

    \node [labeled_vertex, label={[text=red]above:+2}] (4) at (-2, 2) {\(v_1\)};
    \node [labeled_vertex, label={[text=red]above:-6}] (5) at (0.5, 2) {\(v_3\)};
    \node [labeled_vertex, label={[text=red]below:-2}] (7) at (0.5, 0) {\(v_7\)};
    \node [labeled_vertex, label={[text=red]below:+6}] (6) at (-2, 0) {\(v_5\)};

    \node [labeled_vertex, label={[text=red]above:-2}] (8) at (3, 2) {\(v_2\)};
    \node [labeled_vertex, label={[text=red]above:+6}] (9) at (5.5, 2) {\(v_4\)};
    \node [labeled_vertex, label={[text=red]below:+2}] (11) at (5.5, 0) {\(v_8\)};
    \node [labeled_vertex, label={[text=red]below:-6}] (10) at (3, 0) {\(v_6\)};

    \draw[arc] (4) to (5);
    \draw[arc] (5) to (7);
    \draw[arc] (7) to (6);
    \draw[arc] (6) to (4);

    \draw[arc] (10) to (11);
    \draw[arc] (11) to (9);
    \draw[arc] (9) to (8);
    \draw[arc] (8) to (10);
\end{tikzpicture}
\caption{ $2\ar{C_4}$, with  nodes labeled  $2$ and $6$.}
\label{fig:thm1L4b}
\end{subfigure}

\caption{Two labelings of the $2-$node ornament $2\ar{C_4}$, where  $h(v_i)=n+i$, $i\in \{1,2, \hdots, 8\}$; here $n\geq3$. The values of wt$_h$ are marked in red next to each vertex. With either variation of the labeling, given any DDMOG $\ar{G_0}$ with $n$ vertices, $\ar{G_0}\oplus_h(2\ar{C_4})$ is a DDMOG. }
\label{fig:combined_thm1}
\end{figure}

\begin{corollary}
 min $\{ L\in \mathbb{Z}^+ \; | \; 4\ar{C_L} \; \text{is a $4-$nodes ornament} \} = 4$. \end{corollary}  
 
\begin{proof}
Using the two variants of the labeling of $2\ar{C_4}$ from Example~\ref{ex:2nodes2C4}, 
we can construct a $4-$nodes ornament $4\ar{C_L}$ for $L=4$, with nodes labeled $1,2,3,6$. The labels   in one set of $2\ar{C_4}$ will be shifted by $8$ to ensure distinct labels for all the vertices in $4\ar{C_4}$. 
This  lowers the minimal known value for $L$ as per Theorem~\ref{theorem:4CL-ornaments} (as listed in  Table~\ref{s-table}, that value was $L=7$).
  That is,  $\min \{L \; | \; 4\ar{C_L} \; \text{is a   $4-$nodes ornament } \} \leq 4$.
  {The proof is complete if we  show that $4C_3$ cannot be a $4-$nodes ornament. } 
  
  We  assume the opposite: Suppose $4\ar{C_3}$ is a $4-$nodes ornament. Then at the attachment nodes $u_i$, $1\leq i \leq 4$ there should be more than two newly added edges (since otherwise the new weight values at the nodes will not be zeroes). Note that we have only $12$ vertices in $4\ar{C_3}$ {(when $\mathrm{wt}_h(v) = p$, the weighted sum definition only adds new edges between $v$ and node $u_{|p|}$. Thus, } there can be exactly three new edges at each node (either two incoming and one outgoing, or one incoming and two outgoing). For each $i$, the new weights that are added at  $u_i$ must satisfy $a_i+b_i=c_i$ - otherwise the new weight value at $u_i$ will not be zero.  However, it is impossible to have a set of $12$ consecutive integer labels that are $\cup_{i=1}^4 \{a_i, b_i, a_i +b_i\}$. Thus, $L$ cannot be $3$. 

\end{proof}
\section{Searching for sparse DDMOGs}
\label{section:sparse-ddmog-search}
In this section, we investigate what the sparsest nontrivial DDMOG with \(n\) vertices is.
In \cite{FM} DDMOGs with \(n\) vertices were enumerated up to \(n \le 8\). We iterated over all distinct non-trivial DDMOGs for \(n \le 10\)
(see Table \ref{table:DDMOG-counts} for the counts and Algorithm \ref{alg:DDMOGIterator} for the method of enumeration).

Algorithm \ref{alg:DDMOGIterator}   iterates through all DDMOG graphs with \(n\) vertices recursively. The idea is to start with the labeled graph \(nK_1\), and for each vertex \(v = 1, 2, \cdots, n\), find every possible way to connect \(v\) to other vertices in the graph such that the weight of \(v\) is 0.
Specifically, on line \ref{alg:DDMOGIterator:begin-recursion} the recursive search begins with the graph \(nK_1\) and vertex \(1\).
On line \ref{alg:DDMOGIterator:init-count}, in the recursive procedure, a variable is initialized to save the number of nontrivial DDMOGs that can be formed by adding arcs to \(G\) between any of the vertices \(v, v+1, \cdots, n\). Since the procedure begins with \(G = nK_1\) and \(v = 1\), on the first call, this variable will eventually contain the total count of all nontrivial DDMOGs with \(n\) vertices.
After this variable is initialized, we check if we have iterated through all the vertices on line \ref{alg:DDMOGIterator:base-case}. If we have iterated through all the vertices, then \(G\) is a valid DDMOG except \(G\) may be trivial. If \(G\) is nontrivial, then we have found exactly 1 nontrivial DDMOG and update the total count accordingly on line \ref{alg:DDMOGIterator:found-nontrivial-DDMOG}.

If there are still vertices to iterate through i.e.\ \(v < n\), then the weight \(w\) of \(v\) in \(G\) is computed on line \ref{alg:DDMOGIterator:compute-weight}, and we begin a search through all possible ways to connect \(v\) to the vertices \(v + 1, v+2, \cdots, n\). 
Note that if \(v > 1\), then it is redundant to consider the ways to connect \(v\) to the vertices \(1, 2, \cdots, v-1\) since such connections will already have been considered by the recursive procedure for smaller values of \(v\).
Each tuple looped through on line \ref{alg:DDMOGIterator:begin-arc-search} represent potential arcs between \(v\) and the other vertices in the graph. 
For example, if \(v = 1\) and the tuple is \([1,-1,0,\cdots, 0]\), we are representing that vertex \(2\) is into \(v\), vertex \(3\) directed out of \(v\), and vertices \(4, 5,\cdots,n\) not directed into or out of \(v\).
On line \ref{alg:DDMOGIterator:check-DDM} we check if the arc connections result in a weight of \(-w\). If this is the case, then the arcs represented by the array ``neutralize'' the weight \(v\) currently has from the arcs between \(v\) and the vertices \(1, 2, \cdots, v-1\), meaning \(v\) will have a weight of \(0\) when the new arcs in the tuple are added to \(G\).

When such an tuple has been found, the graph \(G\) is copied to a new graph \(G'\) on line \ref{alg:DDMOGIterator:copy-G}, and the arcs represented in the tuple are added to the graph \(G'\). The procedure then recurses using \(G'\) and the next vertex on line \ref{alg:DDMOGIterator:recurse} adding to the current count any new nontrivial DDMOGs found using \(G'\) as the graph in the procedure.
Finally, after all possible arcs in or out of \(v\) are considered, the total count is returned from the procedure on line \ref{alg:DDMOGIterator:return-count}. Note that every nontrivial DDMOG found by the algorithm is available on line \ref{alg:DDMOGIterator:found-nontrivial-DDMOG} for further computations. In Table \ref{table:DDMOG-counts} we record the minimal number of arcs \(MS(n)\) for each value of \(n\), not just the count.

\begin{table}[h!]
\centering
\begin{tabular}{c | c c c c c c c}
     \(n\) & \(<5\) & 5 & 6 & 7 & 8 & 9 & 10 \\
     \hline
     Count & 0 & 6 & 22 & 296 & 9930 & 754804 & 130528634 \\
     \hline
     \(MS(n)\) & Undefined & 8 & 11 & 12 & 13 & 15 & 16
\end{tabular}
\caption{The number of non-trivial DDMOGs with \(n \le 10\) vertices 
}
\label{table:DDMOG-counts}
\end{table}
\subsection{Minimal Sparsity}
As observed in Property 1 in \cite{FM}, every vertex in a DDMOG needs to be incident to at least \(3\) arcs, not all directed into or away from the vertex. The degree sum formula then gives us a lower bound on the number of edges in a DDMOG of \(\ceil{\frac{3n}{2}}\). However, when \(n \equiv 2 \pmod{4}\), we cannot have a \(3\)-regular DDMOG forcing the lower bound to be \(\ceil{\frac{3n}{2}} + 1\). We state these facts in Lemma \ref{lemma:minimum-edges} below.

\begin{lemma}
\label{lemma:minimum-edges}
    Suppose \(\vec{G}\) is a non-trivial DDMOG with \(n\) vertices.
    \[
    \abs{E(\vec{G})} \ge \begin{cases}
    \ceil{\frac{3n}{2}} + 1 & \text{if } n \equiv 2 \pmod{4} \\
    \ceil{\frac{3n}{2}} & \text{otherwise}
    \end{cases}
    \]
\end{lemma}
\begin{proof}
    Property 1 shown in \cite{FM} states that each vertex in a graph with a DDM labeling must have degree at least \(3\). So, if \(G\) has a DDM labeling, then
    \[
    2\abs{E(G)} = \sum_{v \in V(G)} \deg(v) \ge \sum_{v \in V(G)} 3 = 3n.
    \]
    So, \(\abs{E(G)} \ge \ceil{\frac{3n}{2}}\) since the number of edges must be an integer.
    To see that \(\abs{E(G)}\) cannot be \(\frac{3n}{2}\) when \(n \equiv 2 \pmod{4}\), suppose for the sake of contradiction that this does happen.
    Noting that \(\frac{3n}{2}\) is an integer, the degree sum formula implies:
\[
\frac{1}{n}\sum_{v \in V(G)} \deg(v) = \frac{1}{n}\cdot2\abs{E(G)} = \frac{1}{n}\cdot2\left(\frac{3n}{2}\right) = 3.
\]
Since the degree of each vertex must be at least \(3\), it follows that \(G\) must be \(3\)-regular. Let the following system represent the weight equations of each vertex in \(G\).
\[
\begin{cases}
0 = x_1 - x_2 - x_3 \\
0 = x_4 - x_5 - x_6 \\
\cdots \\
0 = x_{3n - 2} - x_{3n - 1} - x_{3n}
\end{cases}
\]
Since for any vertex \(v\), the label of \(v\) must appear
in exactly \(3\) weight equations,
adding up every label we obtain the following equation.
\[
\sum_{i=1}^{3n} x_i = 3 \sum_{v \in V(G)} \ell(v) = 3\frac{n(n+1)}{2}.
\]
Observe:
\[
\sum_{i=1}^{3n} x_i = \sum_{i=1}^n (x_{3i - 2} + x_{3i - 1} + x_{3i}) = \sum_{i=1}^n\left(x_{3i-2} + (x_{3i-2})\right)  = 2\sum_{i=1}^{n} x_{3i-2}.
\]
Setting both expressions of \(\sum_{i=1}^{3n} x_i\) equal, we find
\(\displaystyle
4\sum_{i=1}^{n} x_{3i-2} = 3n(n+1)
\).
Hence \(4\) divides \(3n(n+1)\). However, \(4\) divides \(3n(n+1)\) only when \(n \equiv 0, 3 \pmod{4}\) contradicting that \(n \equiv 2 \pmod{4}\).
\end{proof}

\begin{figure}[h!]
    \centering
    \begin{tikzpicture}
    \tikzset{
        labeled_vertex/.style = {circle, draw=black, fill=white, inner sep = 1.5pt, text=black, minimum size=20pt, very thick},
        arc/.style = {->, very thick}
    };
    \begin{scope}[shift={(-1.75,0)}]
        \node [labeled_vertex] (7) at (0,0) {\(7\)};
        \node [labeled_vertex] (10) at (1,1) {\(10\)};
        \node [labeled_vertex] (1) at (-1,1) {\(1\)};
        \node [labeled_vertex] (3) at (-1,-1) {\(3\)};
        \node [labeled_vertex] (8) at (1,-1) {\(8\)};
        \draw[arc] (7) to (3);
        \draw[arc] (7) to (8);
        \draw[arc] (10) to (7);
        \draw[arc] (1) to (7);
        \draw[arc] (10) to (1);
        \draw[arc] (1) to (3);
        \draw[arc] (3) to (8);
        \draw[arc] (8) to (10);
    \end{scope}
    \begin{scope}[shift={(1.75,0)}]
        \node [labeled_vertex] (4) at (0,0) {\(4\)};
        \node [labeled_vertex] (6) at (1,1) {\(6\)};
        \node [labeled_vertex] (5) at (-1,1) {\(5\)};
        \node [labeled_vertex] (2) at (-1,-1) {\(2\)};
        \node [labeled_vertex] (9) at (1,-1) {\(9\)};
        \draw[arc] (4) to (2);
        \draw[arc] (4) to (9);
        \draw[arc] (6) to (4);
        \draw[arc] (5) to (4);
        \draw[arc] (6) to (5);
        \draw[arc] (5) to (2);
        \draw[arc] (2) to (9);
        \draw[arc] (9) to (6);
    \end{scope}
\end{tikzpicture}
    \caption{\(2W_4\) is a minimal sparsity \(10\) vertex DDMOG}
    \label{fig:order-10-sparse-disconnected-DDMOG}
\end{figure}
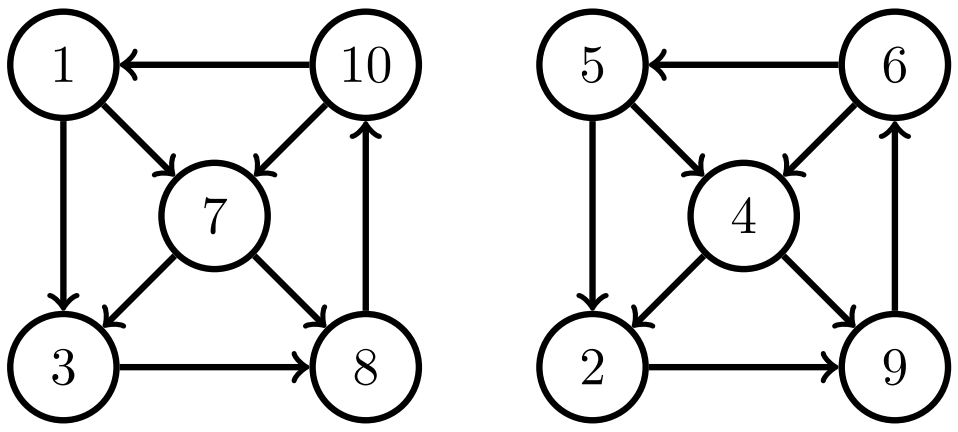

When \(n = 10\) notice that the lower bound in Lemma \ref{lemma:minimum-edges} is reached (Figures \ref{fig:order-10-sparse-disconnected-DDMOG} and \ref{fig:order-10-sparse-connected-DDMOG} illustrate examples).
When \(n \ge 10\), Proposition \ref{proposition:minimal-sparsity-families} below gives examples of nontrivial DDMOGs that also attain this minimum.

\begin{proposition}
\label{proposition:minimal-sparsity-families}
The graphs in table \ref{table:minimal-sparsity-families} are DDMO. Moreover, for each graph \(G\) in the ``Base'' column, the graph \(xK_{3,3} \sqcup G\) is also DDMO for every nonnegative even integer \(x\).
If there are two numbers in the left column (e.g.\ the first row), then \(xK_{3,3} \sqcup G\) is DDMO for every nonnegative integer \(x\) (not just even values).

\begin{table}[h!]
\centering
\begin{tabular}{ c | c | c }
\hline
\(n \bmod 12\) & Base & Figure \\
\hline
1,7 & \(O\) & \ref{fig:order-13-sparse-connected-DDMOG} \\
2 & \(Q\sqcup K_{3,3}\) & \ref{fig:order-14-sparse-disconnected-DDMOG} \\
3,9 & \(H \sqcup W_4\) & \ref{fig:order-15-sparse-disconnected-DDMOG} \\
4 & \( H\sqcup K_{3,3}\) & \ref{fig:order-16-sparse-disconnected-DDMOG} \\
5,11   & \(W_4\) & \ref{fig:a} \\
6 & \(H \sqcup Q\) & \ref{fig:order-18-sparse-disconnected-DDMOG} \\
8 & \(2H \sqcup 2K_{3,3}\) & \ref{fig:order-32-sparse-disconnected-DDMOG} \\
10 & \(2W_4\) & \ref{fig:order-10-sparse-disconnected-DDMOG} \\
\hline
\end{tabular}
\caption{Some minimal sparsity families of DDMOGs. For each graph \(G\) in the ``Base'' column, the graph \(xK_{3,3}\sqcup G\) is DDMO for every even \(x\) (unless there are two numbers in the \(n \bmod 12\) column, then \(xK_{3,3} \sqcup G\) is DDMO for every \(x\))}
\label{table:minimal-sparsity-families}
\end{table}
\end{proposition}

\begin{proof}
    The proof of this Proposition consists of two parts: Theorem \ref{theorem:G2xK_3,3isDDMO} for \(n > 138\) and labelings found via brute force of \(G\sqcup 2x K_{3,3}\) available online at \cite{ZW} for \(n \le 138\). Specifically, if graphs with more than 138 vertices can be constructed by using graphs with \(n = 10, 11, \cdots, 21\) vertices as the graph \(G\) in Theorem \ref{theorem:G2xK_3,3isDDMO}.
    
\end{proof}

Proposition \ref{proposition:minimal-sparsity-families} is the generalization of examples found after loading the definition of a non-trivial DDMOG into a SAT solver (see \cite{ZW} for code implementations)
along with the constraint that the graph attains the minimum in Lemma \ref{lemma:minimum-edges}.

\subsection{A particular orientation of \texorpdfstring{\(K_{3,3}\)}{K₃,₃}}
Many of the example DDMOGs generated via the programs in \cite{ZW} consist of the disjoint union of some base graph along with several copies of \(K_{3,3}\) oriented as illustrated in Figure \ref{fig:oriented_K_3,3}.

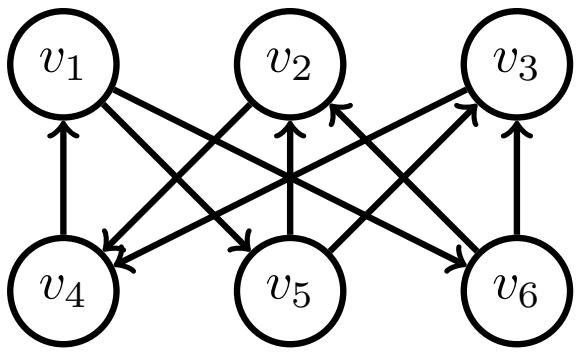
\begin{figure}[h!]
    \centering
    \begin{tikzpicture}
        \tikzset{
        labeled_vertex/.style = {circle, draw=black, fill=white, inner sep = 1.5pt, text=black, minimum size=20pt, very thick},
        arc/.style={->,very thick}
        };
        \node [labeled_vertex] (v4) at (0, 0) {\(v_4\)};
        \node [labeled_vertex] (v5) at (1.5, 0) {\(v_5\)};
        \node [labeled_vertex] (v6) at (3, 0) {\(v_6\)};
        \node [labeled_vertex] (v1) at (0, 1.5) {\(v_1\)};
        \node [labeled_vertex] (v2) at (1.5, 1.5) {\(v_2\)};
        \node [labeled_vertex] (v3) at (3, 1.5) {\(v_3\)};
        \draw[arc] (v4) to (v1);
        \draw[arc] (v1) to (v5);
        \draw[arc] (v1) to (v6);
        \draw[arc] (v2) to (v4);
        \draw[arc] (v5) to (v2);
        \draw[arc] (v6) to (v2);
        \draw[arc] (v3) to (v4);
        \draw[arc] (v5) to (v3);
        \draw[arc] (v6) to (v3);
    \end{tikzpicture}
    \caption{\(K_{3,3}\) with a particular orientation}
    \label{fig:oriented_K_3,3}
\end{figure}

\par\noindent
Letting \(\ell\) be the label function of this graph, note the weight equation of each vertex below.
\[
\begin{cases}
    wt(v_1) = \ell(v_4) - \ell(v_5) - \ell(v_6) \\
    wt(v_2) = -\ell(v_4) + \ell(v_5) + \ell(v_6) \\
    wt(v_3) = -\ell(v_4) + \ell(v_5) + \ell(v_6) \\
    wt(v_4) = -\ell(v_1) + \ell(v_2) + \ell(v_3) \\
    wt(v_5) = \ell(v_1) - \ell(v_2) - \ell(v_3) \\
    wt(v_6) = \ell(v_1) - \ell(v_2) - \ell(v_3)
\end{cases}
\]
If we assume that the weight of each vertex is 0, this system
reduces to the following.
\[
\begin{cases}
    0 = \ell(v_1) - \ell(v_2) - \ell(v_3) \\
    0 = \ell(v_4) - \ell(v_5) - \ell(v_6) 
\end{cases}
\]
So, if \(\ell\) was a DDM labeling,
each integer in \(\{1,2,3,4,5,6\}\) would
be used exactly once in the above system.
Unfortunately, this is not possible, and this agrees
with Theorem \ref{theorem:K_p,q-is-DDMO}.
However, if we consider \(2K_{3,3}\) instead, where
both copies of \(K_{3,3}\) have the orientation
in Figure \ref{fig:oriented_K_3,3}, then if each
weight is 0, the system
of weight equations reduces to the following.
\[
\begin{cases}
    0 = \ell(v_1) - \ell(v_2) - \ell(v_3) \\
    0 = \ell(v_4) - \ell(v_5) - \ell(v_6) \\
    0 = \ell(v_7) - \ell(v_8) - \ell(v_9) \\
    0 = \ell(v_{10}) - \ell(v_{11}) - \ell(v_{12})
\end{cases}
\]
In this system \(v_1, \cdots, v_6\) are the vertices
in one copy of \(K_{3,3}\) and \(v_7, \cdots, v_{12}\)
are the vertices in the other copy of \(K_{3,3}\).
Similar to before, if the labeling is DDM,
then each integer in \(\{1,2,\cdots, 12\}\)
must be used exactly once in this system.

\begin{definition}
\label{definition:difference-triples}
    A difference triple is a tuple of 3 positive integers \((a,b,c)\) such that \(a + b = c\).
\end{definition}

A solution of the previous system of equations then corresponds to partitioning the array \((1,2,\cdots, 12\) into difference triples.
Such a partition can be constructed using Langford sequences.

\begin{lemma}
    \label{lemma:langford_triples}
If an order \(n\) Langford sequence with defect \(d\) exists,
then \(A = (d, d + 1, \cdots, d + 3n - 1)\)
can be partitioned into difference triples.
\end{lemma}

\begin{proof}
    Let \(\{(u_i, v_i)\}_{i=1}^k\) be an order \(n\) Langford pairing with defect \(d\),
    and define \(a_i = v_i - u_i\), \(b_i = u_i + n + d - 1\), and \(c_i = a_i + b_i\).

    It remains to show that \(A = \bigcup_{i=1}^n \{a_i, b_i, c_i\}\). To see this,
    recall from the definition of a Langford pairing
    that \(d \le a_i \le n + d - 1\). Observe:
    \[
    \bigcup_{i=1}^n\{\underbrace{u_i + n + d - 1}_{b_i}, \underbrace{v_i + n + d - 1}_{c_i}\} = \{n + d, n + d + 1,\cdots, 3n + d - 1\}.
    \]
\end{proof}

\begin{remark}
    The construction used in Lemma \ref{lemma:langford_triples} is a minor variation of the one first conceived by Skolem in \cite{skolemRemarksTripleSystems1958} where an order \(n\) Skolem sequence was used to create Steiner triple systems by partitioning the array \((1,2,\cdots, 3n)\) into difference triples.
\end{remark}

\begin{theorem}
    \label{theorem:G2xK_3,3isDDMO}
Let \(G\) be a DDMO graph with \(n\) vertices and \(x\) be some non-negative integer.
If \(2x > n\), then \(G \sqcup 2xK_{3,3}\) is DDMO.
\end{theorem}
\begin{proof}
Lemma \ref{lemma:langford-sequence-existence} guarantees that a Langford sequence of order \(4x\) exists whenever \(4x \ge 2d - 1\).
By choosing the defect to be \(n + 1\), then
this inequality becomes \(4x \ge 2n + 1\) or equivalently:
\(2x > n\).
Lemma \ref{lemma:langford_triples} guarantees that we can then partition the sequence
\((n+1, n+2, \cdots, n+1 + 12x)\) into difference triples.
Using the DDMO labeling of \(G\), and by assigning each difference triple to an instance of \(K_{3,3}\) (as was done in Figure \ref{fig:oriented_K_3,3}),
we obtain a DDMO
labeling of \(G\sqcup 2xK_{3,3}\).
\end{proof}

By using minimally sparse DDMOGs as our base graph \(G\), we can construct arbitrarily large minimally sparse DDMOGS.
Examples of minimally sparse DDMOGs for \(n = 10, \cdots, 21\) exist according to the case of \(n < 138\) in Proposition \ref{proposition:minimal-sparsity-families}.

\section{Conclusion}\label{sec:final}

We have shown various constructions for ornaments that are focused around the ideas of minimality and sparsity. To this, we offer 3 theorems and various lemmas, propositions, and corollaries aided by examples to convey the new findings and how they differ from what has been found previously. In Section 2, we discussed different configurations and labelings of ornaments to reduce the number of cycles needed to establish a DDMOG. Then in Section 3, we brought the idea of sparsity in DDMOGs and ornaments to the forefront and discussed an algorithm to do so, while focusing on particular $K_{3,3}$ graphs and Langford sequences. 

The ornaments in this paper are attached via the weighted sum construction. In Definition \ref{definition:general-ornament} below we generalize this to any attachment scheme.
\begin{definition}
    \label{definition:general-ornament}
    Let \(\ar{G}\) be a DDMOG with a DDM labeling \(f_0: V(\ar{G}) \rightarrow \{1,\cdots,|{V(\ar{G})}|\}\).
    An oriented graph \(\ar{H}\) is an \(s\)-node ornament attachable to \(\ar{G}\) if and only if there exists a bijective labeling \(h: V(\ar{H}) \rightarrow \{n+1,\cdots, n+|{V(\ar{H})}|\}\) 
    such that the function \(f: V(\ar{G})\cup V(\ar{H}) \rightarrow \{1,\cdots,|{V(\ar{G})}| + |{V(\ar{H})}|\}\) given by \begin{equation}\label{deffonG0H} f(v)=\begin{cases}
			h(v), & \text{if $v \in V(H)$,}\\
            f_0(v), & \text{if $v \in V(G_0)$.}
		 \end{cases}\end{equation}
    is a DDM labeling of the graph \(\ar{G}\cup_E\ar{H}\) (the graph with vertices \(V(\ar{G}) \cup V(\ar{H})\) and arcs \(E(\ar{G}) \cup E(\ar{H}) \cup E\)).
\end{definition}

In \cite{ALR}, they asked what is the sparsest \(s\)-node ornament with \(k\) vertices is.
A graph with \(3K_1\) attached as in Proposition \ref{proposition:3mK_1-ornament} below does not have minimal sparsity. However, the generalized ornament itself \(3K_1\) has zero edges.

In Figure \ref{fig:oriented_K_3,3}, we connected two sets of 3 vertices together in such a way where if each part of \(K_{3,3}\) is labeled using a difference triple, then the weight of each vertex is 0. This connection strategy can be used to attach generalized ornaments. The following proposition gives an example of how to make \(3mK_1\) an ornament.

\begin{proposition}
    \label{proposition:3mK_1-ornament}
    Let \(\ar{G}\) be a DDMOG with \(n \ge 3\) vertices.
    If an order \(m\) Langford sequence exists with defect \(n+1\), then \(3mK_1\) is an \(s\)-node generalized ornament attachable to \(\ar{G}\) for each \(s = 3, 4, \cdots, n\). 
\end{proposition}
\begin{proof}
Suppose that a Langford sequence of order \(m\) exists with defect \(n+1\) and label \(3mK_1\) with the integers \(n+1, n+2, \cdots, n +3m\).
Proposition \ref{lemma:langford_triples}
then guarantees that \((n+1, n+2,\cdots, n+3m)\) can be partitioned into difference triples \(\{(a_i, b_i, c_i)\}_{i=1}^m\).

Notice that for each \(k \ge 3\) that the tuple \((1, k-1, k)\) is a difference triple. So, if \(3 \le s \le n\), then \(\{(1, k-1, k)\}_{k=3}^s\) is a collection of difference triples formed solely from the labels in \(\ar{G}\).
Since Lemma \ref{lemma:langford-sequence-existence} forces \(m \ge 2(n+1) - 1 \ge n\), each set of 3 vertices in \(\ar{G}\) corresponding to a difference triple \((1, k-1,k)\) can be connected to at least one set of 3 vertices in \(3mK_1\) corresponding to a difference triple \((a_i, b_i, c_i)\) with arcs like in Figure \ref{fig:oriented_K_3,3}.
\end{proof}

\section*{Code Availability}
The implementation of Algorithm \ref{alg:DDMOGIterator} used in this article is publicly available in the GitHub repository at \url{https://github.com/zachwilcher/DDMOG/blob/main/ddmog_counter/main.cpp}
 
 \section*{Acknowledgments}

This work utilized the Ball State University beowulf cluster, which is supported
by The National Science Foundation (MRI-1726017) and Ball State University,
Muncie, Indiana. M.Lort was supported by a Pepsi Undergraduate Research Grant from Ball State University.

We thank Dr.\ Justin Murray and Dr.\ John Lorch for their helpful discussions.  The first author thanks the Institute for Computational and Experimental Research in Mathematics (ICERM) for the invitation to participate  in its Research Experiences for Undergraduate Faculty (REUF) program for two consecutive years.
\vspace{5mm}


\newpage
\section*{Appendix}

  \begin{figure}[htbp]
\centering

\begin{subfigure}[b]{\textwidth}
  \centering
  \begin{tikzpicture}[scale=0.76]
    \tikzset{
        labeled_vertex/.style = {circle, draw=black, fill=white, inner sep = 1.5pt, text=black, minimum size=20pt, very thick},
        arc/.style={->,very thick}
    };

    \node [labeled_vertex] (1) at (1, 6.5) {1}; 

    \node [labeled_vertex] (5) at (-1, 5) {5};
    \node [labeled_vertex] (3) at (3, 5) {3};
    \node [labeled_vertex] (2) at (-1, 8) {2};
    \node [labeled_vertex] (4) at (3, 8) {4};

    \draw[arc] (5) to (1);
    \draw[arc] (5) to (2);
    \draw[arc] (3) to (5);

    \draw[arc] (2) to (4);
    \draw[arc] (4) to (3);
    \draw[arc] (2) to (1);
    \draw[arc] (1) to (3);
    \draw[arc] (1) to (4);
  \end{tikzpicture}
  
  \vspace{0.2cm}
  \caption{DDM Labeling of $\overrightarrow{W_4}$ ($\overrightarrow{W_4}=\overrightarrow{K_1}\oplus \overrightarrow{C_4}$)}
  \label{fig:a}
\end{subfigure}

\vspace{1cm} 

\begin{subfigure}[b]{\textwidth}
  \centering
  \begin{tikzpicture}[scale=0.76]
    \tikzset{
        labeled_vertex/.style = {circle, draw=black, fill=white, inner sep = 1.5pt, text=black, minimum size=20pt, very thick},
        arc/.style={->,very thick}
    };

    \node [labeled_vertex, label={[text=teal]above left:+1}] (6) at (-11, 2) {6};
    \node [labeled_vertex, label={[text=red]left:-3}] (7) at (-12, 0.5) {7};
    \node [labeled_vertex, label={[text=orange]below left:-4}] (9) at (-11, -1) {9};
    \node [labeled_vertex, label={[text=teal]below right:-1}] (11) at (-9, -1) {11};
    \node [labeled_vertex, label={[text=red]above right:+3}] (10) at (-8, 0.5) {10};
    \node [labeled_vertex, label={[text=orange]above right:+4}] (8) at (-9, 2) {8};

    \draw[arc] (6) to (7);
    \draw[arc] (7) to (9);
    \draw[arc] (9) to (11);
    \draw[arc] (11) to (10);
    \draw[arc] (10) to (8);
    \draw[arc] (8) to (6);

    \node [labeled_vertex, label={[text=teal]above left:-1}] (12) at (-3, 2) {12};
    \node [labeled_vertex, label={[text=red]left:+3}] (13) at (-4, 0.5) {13};
    \node [labeled_vertex, label={[text=orange]below left:+4}] (15) at (-3, -1) {15};
    \node [labeled_vertex, label={[text=orange]above right:-4}] (14) at (-1, 2) {14};
    \node [labeled_vertex, label={[text=red]right:-3}] (16) at (0, 0.5) {16};
    \node [labeled_vertex, label={[text=teal]below right:+1}] (17) at (-1, -1) {17};

    \draw[arc] (12) to (14);
    \draw[arc] (14) to (16);
    \draw[arc] (16) to (17);
    \draw[arc] (17) to (15);
    \draw[arc] (15) to (13);
    \draw[arc] (13) to (12);
  \end{tikzpicture}
  
  \vspace{0.2cm}
  \caption{Labeling of the $3$-nodes ornament $2\overrightarrow{C_6}$ (values of $\text{wt}_h$: \textcolor{teal}{$\pm 1$}, \textcolor{red}{$\pm 3$} and \textcolor{orange}{$\pm 4$})}
  \label{fig:b}
\end{subfigure}

\vspace{1cm} 

\begin{subfigure}[b]{\textwidth}
  \centering
  \begin{tikzpicture}[scale=0.76]
    \tikzset{
        labeled_vertex/.style = {circle, draw=black, fill=white, inner sep = 1.5pt, text=black, minimum size=20pt, very thick},
        arc/.style={->,very thick}
    };

    \node [labeled_vertex] (1) at (0, 6.5) {1}; 
    \node [labeled_vertex] (5) at (-2, 5) {5};
    \node [labeled_vertex] (3) at (2, 5) {3};
    \node [labeled_vertex] (2) at (-2, 8) {2};
    \node [labeled_vertex] (4) at (2, 8) {4};

    \draw[arc] (5) to (1);
    \draw[arc] (5) to (2);
    \draw[arc] (3) to (5);
    \draw[arc] (2) to (4);
    \draw[arc] (4) to (3);
    \draw[arc] (2) to (1);
    \draw[arc] (1) to (3);
    \draw[arc] (1) to (4);

    \node [labeled_vertex] (6) at (-4, 0) {6};
    \node [labeled_vertex] (7) at (-5, -1.5) {7};
    \node [labeled_vertex] (9) at (-4, -3) {9};
    \node [labeled_vertex] (11) at (-2, -3) {11};
    \node [labeled_vertex] (10) at (-1, -1.5) {10};
    \node [labeled_vertex] (8) at (-2, 0) {8};

    \draw[arc] (6) to (7);
    \draw[arc] (7) to (9);
    \draw[arc] (9) to (11);
    \draw[arc] (11) to (10);
    \draw[arc] (10) to (8);
    \draw[arc] (8) to (6);

    \node [labeled_vertex] (12) at (-11, 4) {12};
    \node [labeled_vertex] (13) at (-12, 2.5) {13};
    \node [labeled_vertex] (15) at (-11, 1) {15};
    \node [labeled_vertex] (14) at (-9, 4) {14};
    \node [labeled_vertex] (16) at (-8, 2.5) {16};
    \node [labeled_vertex] (17) at (-9, 1) {17};

    \draw[arc] (12) to (14);
    \draw[arc] (14) to (16);
    \draw[arc] (16) to (17);
    \draw[arc] (17) to (15);
    \draw[arc] (15) to (13);
    \draw[arc] (13) to (12);

    \draw[arc, teal] (6) to [bend left=3] (1);
    \draw[arc, red] (3) to [bend left=39] (7);
    \draw[arc, orange] (8) to (4);
    \draw[arc, teal] (1) to [bend left=2] (11);
    \draw[arc, red] (10) to [bend right=20] (3);
    \draw[arc, orange] (4) to [bend right=7] (9);

    \draw[arc, teal] (1) to [bend right=23] (12);
    \draw[arc, red] (13) to (3);
    \draw[arc, orange] (15) to [bend left=13] (4);
    \draw[arc, teal] (17) to [bend right=49] (1);
    \draw[arc, red] (3) to (16);
    \draw[arc, orange] (4) to [bend right=9] (14);
  \end{tikzpicture}
  
  \vspace{0.2cm}
  \caption{DDM Labeling of $\overrightarrow{W_4} \oplus 2\overrightarrow{C_6}$; the new edges are colored   (corresponding to the values of $\text{wt}_h$)}
  \label{fig:c}
\end{subfigure}

\vspace{0.5cm}
\caption{The construction of DDMOG $\overrightarrow{W_4} \oplus 2\overrightarrow{C_6}$.}
\label{fig:Two6Cycles}
\end{figure}
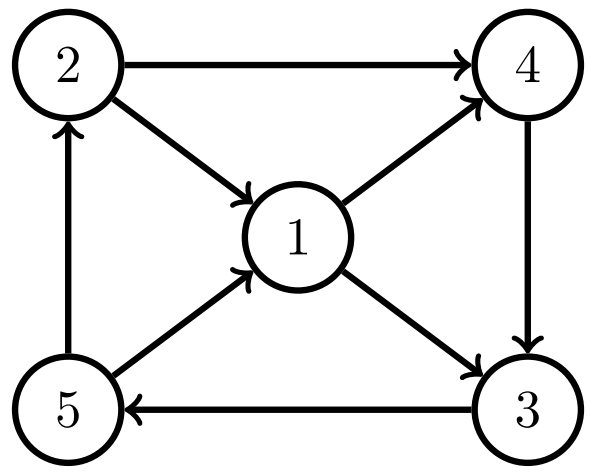
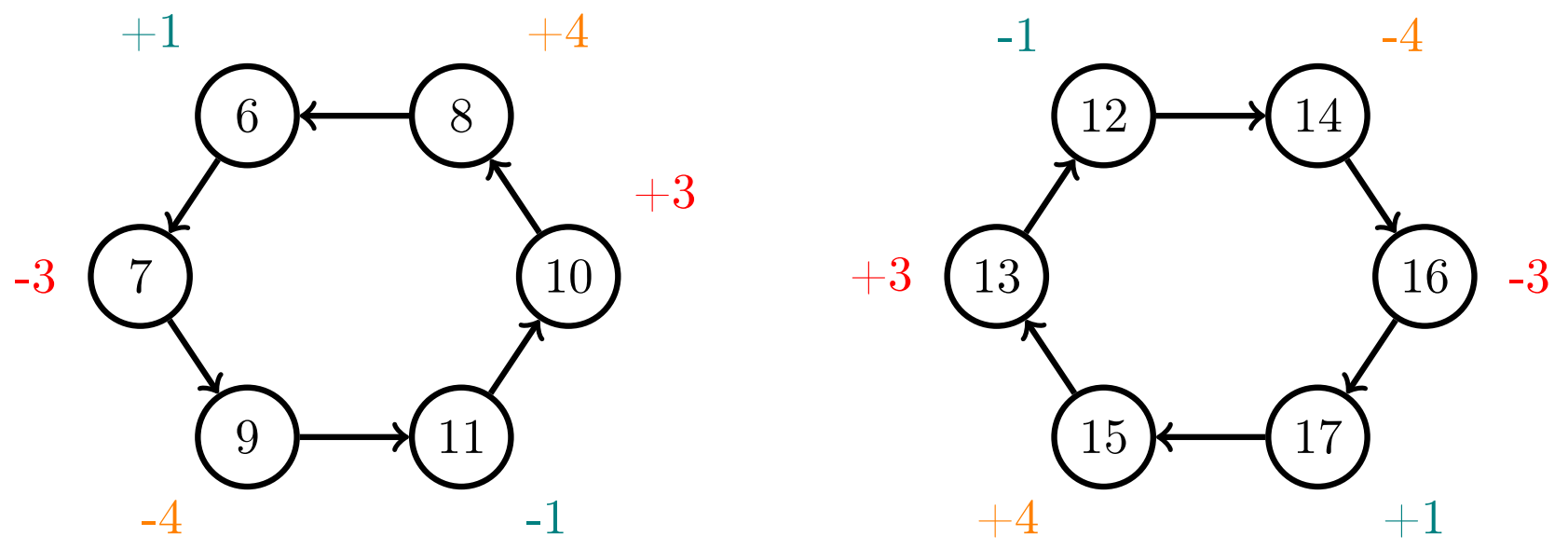
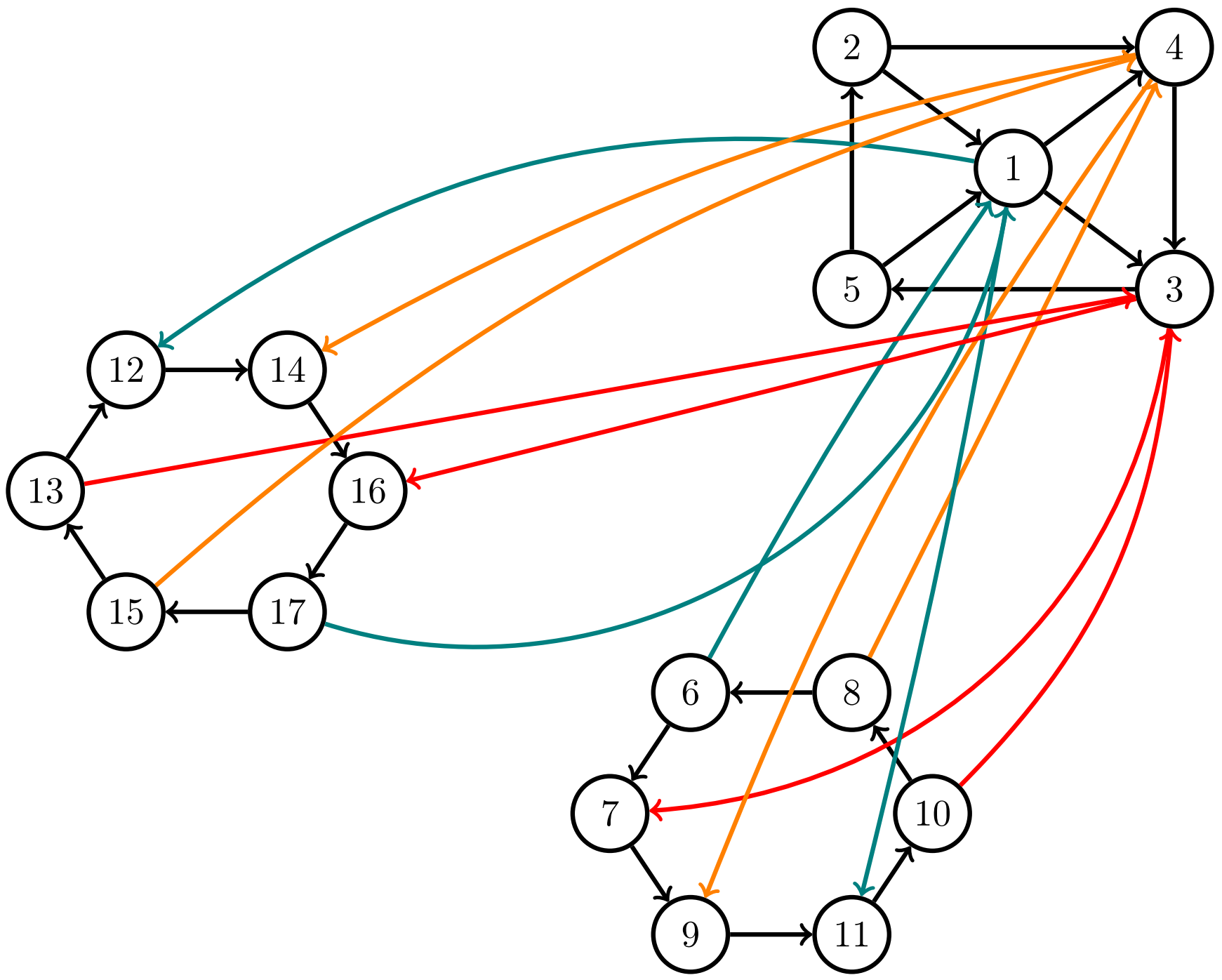

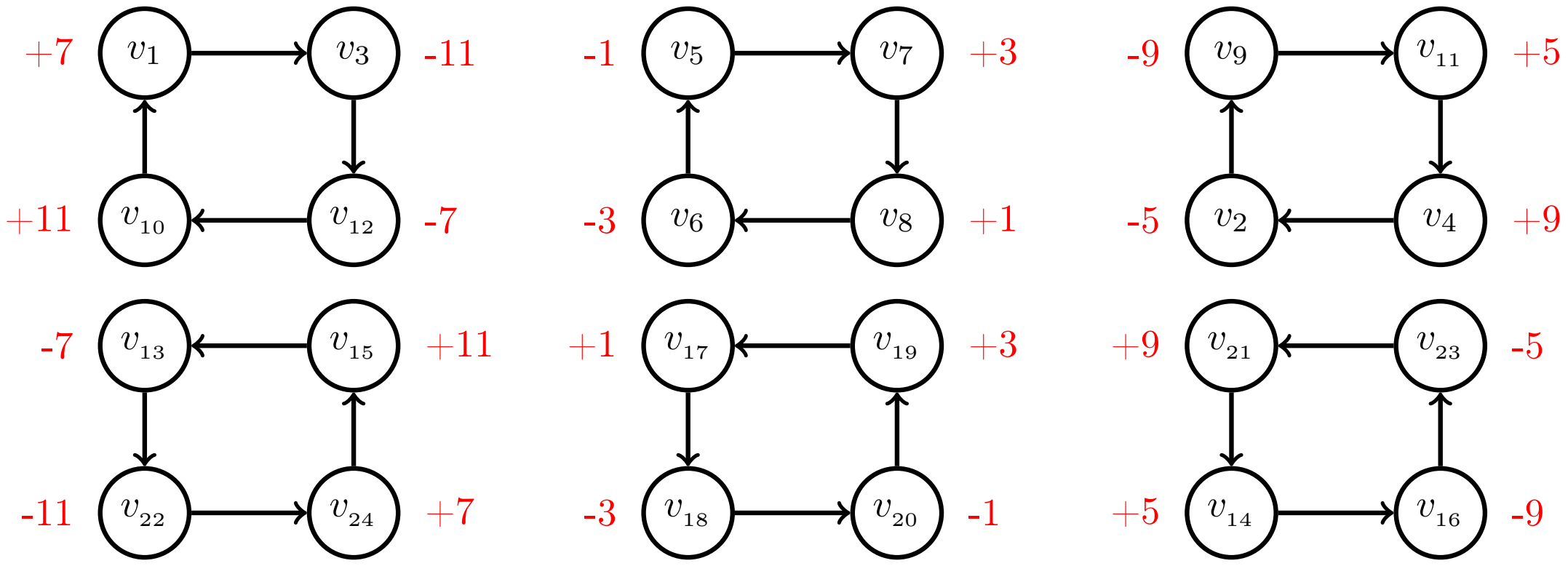
\begin{figure}
\centering
 
 \begin{tikzpicture}[scale=0.76]
    \tikzset{
        labeled_vertex/.style = {circle, draw=black, fill=white, inner sep = 1.5pt, text=black, minimum size=23pt, very thick},
        arc/.style={->,very thick}
    };

    \node [labeled_vertex, label={[text=red]left:+7}] (v1) at (-2.5, 2) {\(v_1\)};
    \node [labeled_vertex, label={[text=red]right:-11}] (v3) at (0, 2) {\(v_3\)}; 
    \node [labeled_vertex, label={[text=red]right:-7}] (v12) at (0, 0) {\(v_{_{12}}\)};
    \node [labeled_vertex, label={[text=red]left:+11}] (v10) at (-2.5, 0) {\(v_{_{10}}\)};
 
    \draw[arc] (v1) to (v3);
    \draw[arc] (v3) to (v12);
    \draw[arc] (v12) to (v10);
    \draw[arc] (v10) to (v1);

    \node [labeled_vertex, label={[text=red]left:-1}] (v5) at (4, 2) {\(v_5\)};
    \node [labeled_vertex, label={[text=red]right:+3}] (v7) at (6.5, 2) {\(v_7\)};
    \node [labeled_vertex, label={[text=red]right:+1}] (v8) at (6.5, 0) {\(v_8\)};
    \node [labeled_vertex, label={[text=red]left:-3}] (v6) at (4, 0) {\(v_6\)};
 
    \draw[arc] (v5) to (v7);
    \draw[arc] (v7) to (v8);
    \draw[arc] (v8) to (v6);
    \draw[arc] (v6) to (v5);

    \node [labeled_vertex, label={[text=red]left:-9}] (v9) at (10.5, 2) {\(v_9\)};
    \node [labeled_vertex, label={[text=red]right:+5}] (v11) at (13, 2) {\(v_{_{11}}\)};
    \node [labeled_vertex, label={[text=red]right:+9}] (v4) at (13, 0) {\(v_4\)};
    \node [labeled_vertex, label={[text=red]left:-5}] (v2) at (10.5, 0) {\(v_2\)};

    \draw[arc] (v9) to (v11);
    \draw[arc] (v11) to (v4);
    \draw[arc] (v4) to (v2);
    \draw[arc] (v2) to (v9);

    \node [labeled_vertex, label={[text=red]left:-7}] (v13) at (-2.5, -1.5) {\(v_{_{13}}\)};
    \node [labeled_vertex, label={[text=red]right:+11}] (v15) at (0, -1.5) {\(v_{_{15}}\)}; 
    \node [labeled_vertex, label={[text=red]right:+7}] (v24) at (0, -3.5) {\(v_{_{24}}\)};
    \node [labeled_vertex, label={[text=red]left:-11}] (v22) at (-2.5, -3.5) {\(v_{_{22}}\)};
 
    \draw[arc] (v13) to (v22);
    \draw[arc] (v22) to (v24);
    \draw[arc] (v24) to (v15);
    \draw[arc] (v15) to (v13);

    \node [labeled_vertex, label={[text=red]left:+1}] (v17) at (4, -1.5) {\(v_{_{17}}\)};
    \node [labeled_vertex, label={[text=red]right:+3}] (v19) at (6.5, -1.5) {\(v_{_{19}}\)};
    \node [labeled_vertex, label={[text=red]right:-1}] (v20) at (6.5, -3.5) {\(v_{_{20}}\)};
    \node [labeled_vertex, label={[text=red]left:-3}] (v18) at (4, -3.5) {\(v_{_{18}}\)};
 
    \draw[arc] (v19) to (v17);
    \draw[arc] (v17) to (v18);
    \draw[arc] (v18) to (v20);
    \draw[arc] (v20) to (v19);

    \node [labeled_vertex, label={[text=red]left:+9}] (v21) at (10.5, -1.5) {\(v_{_{21}}\)};
    \node [labeled_vertex, label={[text=red]right:-5}] (v23) at (13, -1.5) {\(v_{_{23}}\)};
    \node [labeled_vertex, label={[text=red]right:-9}] (v16) at (13, -3.5) {\(v_{_{16}}\)};
    \node [labeled_vertex, label={[text=red]left:+5}] (v14) at (10.5, -3.5) {\(v_{_{14}}\)};

    \draw[arc] (v23) to (v21);
    \draw[arc] (v21) to (v14);
    \draw[arc] (v14) to (v16);
    \draw[arc] (v16) to (v23);

\end{tikzpicture}

\caption{$6-$node ornament, $6\ar{C_4}$; the  labeling function is $h(v_j)=n+j$, $j\in \{1,2,\hdots, 24\}$, given any $n\geq 11$. The values of $\text{wt}_h$ are marked in red next to every vertex.   }
\label{fig:6nodes4vertex}
\end{figure} 

\begin{figure}
    \centering
    \begin{tikzpicture}[scale=3, rotate=33]
   \tikzset{
       labeled_vertex/.style = {circle, draw=black, fill=white, inner sep = 1.5pt, text=black, minimum size=20pt, very thick},
       arc/.style = {->, very thick}
   };
  \node[labeled_vertex] (1) at (-0.60, -0.32) {1};
  \node[labeled_vertex] (2) at (-0.03, -0.68) {2};
  \node[labeled_vertex] (3) at (0.65, 0.38) {3};
  \node[labeled_vertex] (4) at (0.07, 0.75) {4};
  \node[labeled_vertex] (5) at (-0.07, -0.11) {5};
  \node[labeled_vertex] (6) at (-0.60, 0.78) {6};
  \node[labeled_vertex] (7) at (-0.61, 0.29) {7};
  \node[labeled_vertex] (8) at (0.96, -0.21) {8};
  \node[labeled_vertex] (9) at (-0.23, 0.29) {9};
  \node[labeled_vertex] (10) at (0.52, -0.44) {10};
  \node[labeled_vertex] (11) at (-1.00, 0.18) {11};
  \node[labeled_vertex] (12) at (0.35, -0.08) {12};
  \node[labeled_vertex] (13) at (0.59, -0.83) {13};
  \draw[arc] (1) to (11);
  \draw[arc] (2) to (1);
  \draw[arc] (2) to (12);
  \draw[arc] (3) to (12);
  \draw[arc] (4) to (3);
  \draw[arc] (4) to (6);
  \draw[arc] (5) to (9);
  \draw[arc] (5) to (10);
  \draw[arc] (6) to (11);
  \draw[arc] (7) to (5);
  \draw[arc] (7) to (6);
  \draw[arc] (8) to (3);
  \draw[arc] (8) to (10);
  \draw[arc] (9) to (1);
  \draw[arc] (9) to (4);
  \draw[arc] (10) to (13);
  \draw[arc] (11) to (7);
  \draw[arc] (12) to (5);
  \draw[arc] (13) to (2);
  \draw[arc] (13) to (8);
\end{tikzpicture}
    \caption{The \(O\) graph}
    \label{fig:order-13-sparse-connected-DDMOG}
\end{figure}
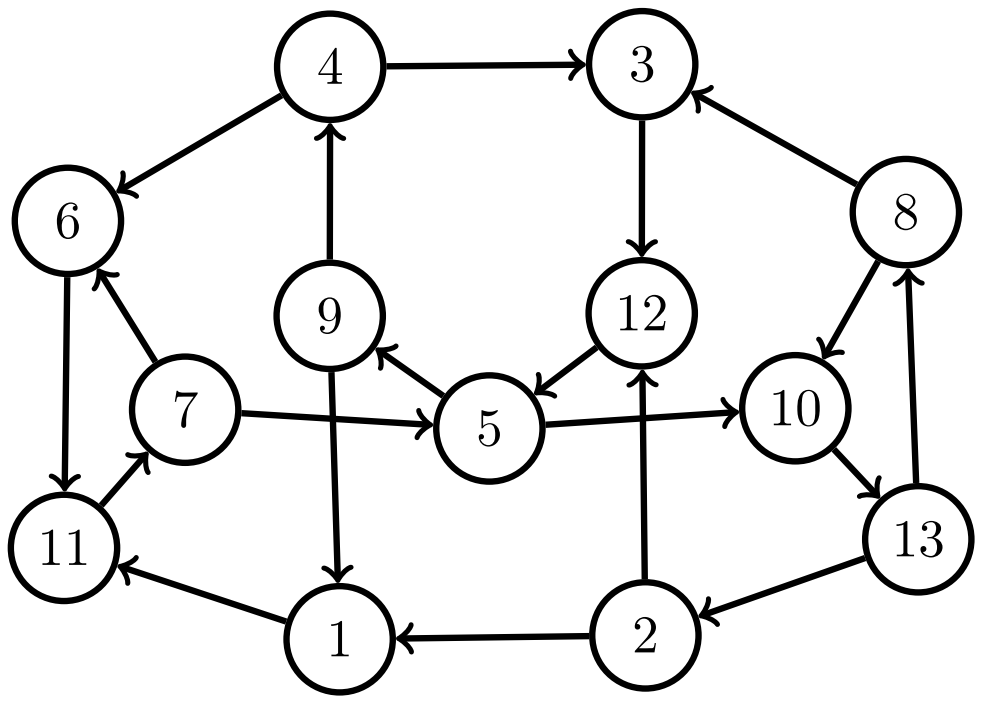

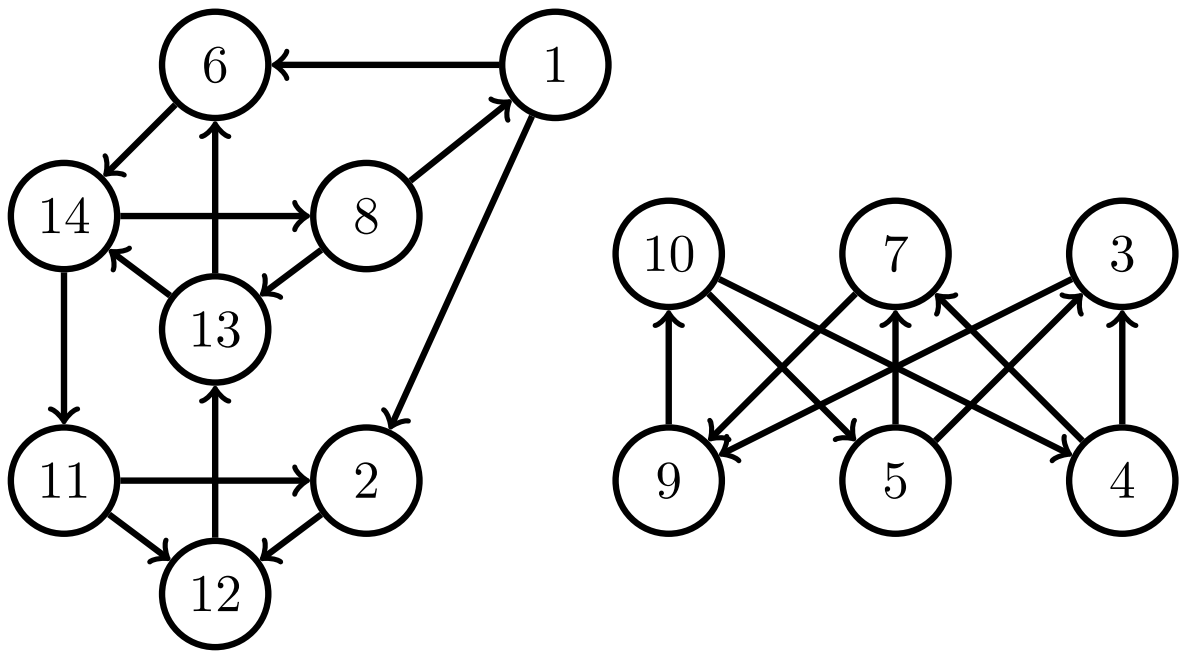
\begin{figure}
    \centering
    \begin{tikzpicture}
    \tikzset{
        labeled_vertex/.style = {circle, draw=black, fill=white, inner sep = 1.5pt, text=black, minimum size=20pt, very thick},
        arc/.style = {->, very thick}
    };
    \begin{scope}[shift={(1,0)}]
        \node [labeled_vertex] (v4) at (0, 0) {\(9\)};
        \node [labeled_vertex] (v5) at (1.5, 0) {\(5\)};
        \node [labeled_vertex] (v6) at (3, 0) {\(4\)};
        \node [labeled_vertex] (v1) at (0, 1.5) {\(10\)};
        \node [labeled_vertex] (v2) at (1.5, 1.5) {\(7\)};
        \node [labeled_vertex] (v3) at (3, 1.5) {\(3\)};
        \draw[arc] (v4) -- (v1);
        \draw[arc] (v1) -- (v5);
        \draw[arc] (v1) -- (v6);
        \draw[arc] (v2) -- (v4);
        \draw[arc] (v5) -- (v2);
        \draw[arc] (v6) -- (v2);
        \draw[arc] (v3) -- (v4);
        \draw[arc] (v5) -- (v3);
        \draw[arc] (v6) -- (v3);
    \end{scope}
    \begin{scope}[shift={(-2,0)}]
        \node [labeled_vertex] (q11) at (-1, 0) {\(11\)};
        \node [labeled_vertex] (q2) at (1, 0) {\(2\)};
        \node [labeled_vertex] (q12) at (0, -.75) {\(12\)};
        \draw[arc] (q11) to (q2);
        \draw[arc] (q2) to (q12);
        \draw[arc] (q11) to (q12);
        \node [labeled_vertex] (q14) at (-1, 1.75) {\(14\)};
        \node [labeled_vertex] (q8) at (1, 1.75) {\(8\)};
        \node [labeled_vertex] (q13) at (0, 1) {\(13\)};
        \draw[arc] (q14) to (q8);
        \draw[arc] (q8) to (q13);
        \draw[arc] (q13) to (q14);
        \draw[arc] (q14) to (q11);
        \draw[arc] (q12) to (q13);
        \node [labeled_vertex] (q6) at (0, 2.75) {\(6\)};
        \draw[arc] (q6) to (q14);
        \draw[arc] (q13) to (q6);
        \node [labeled_vertex] (q1) at (2.25, 2.75) {\(1\)};
        \draw[arc] (q1) to (q6);
        \draw[arc] (q8) to (q1);
        \draw[arc] (q1) to (q2);
        \end{scope}
\end{tikzpicture}
    \caption{\(Q \sqcup K_{3,3}\)}
    \label{fig:order-14-sparse-disconnected-DDMOG}
\end{figure}

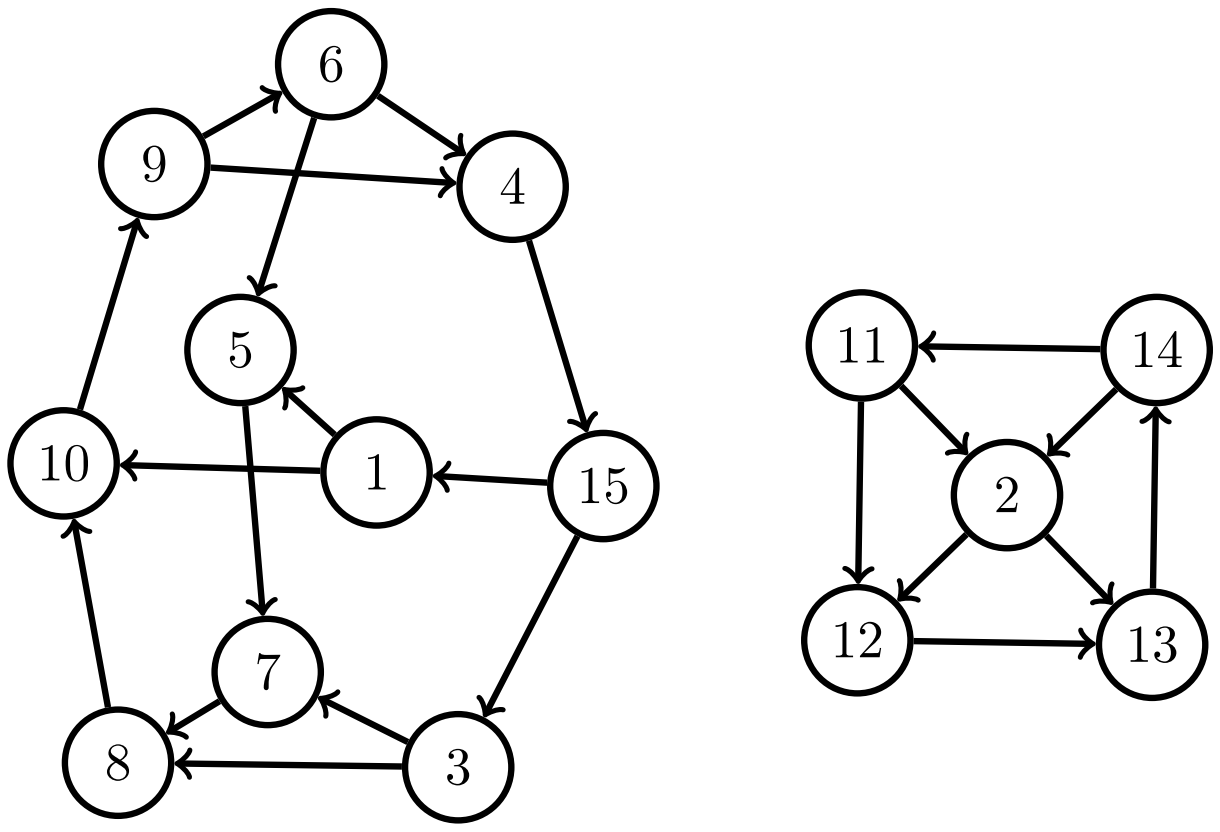
\begin{figure}
    \centering
    \begin{tikzpicture}[scale=3]
  \tikzset{
    labeled_vertex/.style = {circle, draw=black, fill=white, inner sep = 1.5pt, text=black, minimum size=20pt, very thick},
    arc/.style = {->, very thick}
  };
  \begin{scope}[shift={(-.66,0)}]
    \node[labeled_vertex] (1) at (0.07, 0.10) {1};
    \node[labeled_vertex] (3) at (0.25, -0.55) {3};
    \node[labeled_vertex] (4) at (0.37, 0.73) {4};
    \node[labeled_vertex] (5) at (-0.23, 0.37) {5};
    \node[labeled_vertex] (6) at (-0.03, 1.00) {6};
    \node[labeled_vertex] (7) at (-0.17, -0.34) {7};
    \node[labeled_vertex] (8) at (-0.50, -0.54) {8};
    \node[labeled_vertex] (9) at (-0.42, 0.78) {9};
    \node[labeled_vertex] (10) at (-0.62, 0.12) {10};
    \node[labeled_vertex] (15) at (0.57, 0.07) {15};
  \end{scope}

  \begin{scope}[shift={(0.66,0.4)}]
    \node[labeled_vertex] (2) at (0.14, -0.35) {2};
    \node[labeled_vertex] (11) at (-0.18, -0.02) {11};
    \node[labeled_vertex] (12) at (-0.19, -0.67) {12};
    \node[labeled_vertex] (13) at (0.46, -0.68) {13};
    \node[labeled_vertex] (14) at (0.47, -0.03) {14};
  \end{scope}
  \draw[arc] (2) to (12);
  \draw[arc] (2) to (13);
  \draw[arc] (14) to (2);
  \draw[arc] (11) to (2);
  \draw[arc] (14) to (11);
  \draw[arc] (11) to (12);
  \draw[arc] (12) to (13);
  \draw[arc] (13) to (14);
    
  \draw[arc] (1) to (5);
  \draw[arc] (1) to (10);
  \draw[arc] (3) to (7);
  \draw[arc] (3) to (8);
  \draw[arc] (4) to (15);
  \draw[arc] (5) to (7);
  \draw[arc] (6) to (4);
  \draw[arc] (6) to (5);
  \draw[arc] (7) to (8);
  \draw[arc] (8) to (10);
  \draw[arc] (9) to (4);
  \draw[arc] (9) to (6);
  \draw[arc] (10) to (9);
  \draw[arc] (15) to (1);
  \draw[arc] (15) to (3);
\end{tikzpicture}
    \caption{\(H\sqcup W_4\)}
    \label{fig:order-15-sparse-disconnected-DDMOG}
\end{figure}

\begin{figure}
    \centering
    \begin{tikzpicture}[scale=3]
   \tikzset{
       labeled_vertex/.style = {circle, draw=black, fill=white, inner sep = 1.5pt, text=black, minimum size=20pt, very thick},
       arc/.style = {->, very thick}
   };
  \begin{scope}[shift={(0,0)}]
    \node[labeled_vertex] (15) at (0.54, 0.10) {15};
    \node[labeled_vertex] (1) at (0.08, 0.16) {1};
    \node[labeled_vertex] (8) at (0.40, 0.72) {8};
    \node[labeled_vertex] (2) at (-0.18, 0.43) {2};
    \node[labeled_vertex] (5) at (0.04, 1.00) {5};
    \node[labeled_vertex] (9) at (-0.49, -0.40) {9};
    \node[labeled_vertex] (10) at (-0.33, 0.82) {10};
    \node[labeled_vertex] (13) at (-0.55, 0.22) {13};
    \node[labeled_vertex] (6) at (-0.17, -0.23) {6};
    \node[labeled_vertex] (7) at (0.20, -0.45) {7};
  \end{scope}

  \begin{scope}[shift={(1.2,0)},scale=0.33]
    \node [labeled_vertex] (v4) at (0, 0) {\(14\)};
    \node [labeled_vertex] (v5) at (1.5, 0) {\(3\)};
    \node [labeled_vertex] (v6) at (3, 0) {\(11\)};
    \node [labeled_vertex] (v1) at (0, 1.5) {\(12\)};
    \node [labeled_vertex] (v2) at (1.5, 1.5) {\(4\)};
    \node [labeled_vertex] (v3) at (3, 1.5) {\(16\)};
    \draw[arc] (v4) to (v1);
    \draw[arc] (v1) to (v5);
    \draw[arc] (v1) to (v6);
    \draw[arc] (v2) to (v4);
    \draw[arc] (v5) to (v2);
    \draw[arc] (v6) to (v2);
    \draw[arc] (v3) to (v4);
    \draw[arc] (v5) to (v3);
    \draw[arc] (v6) to (v3);
  \end{scope}
  
  \draw[arc] (15) to (1);
  \draw[arc] (15) to (7);
  \draw[arc] (1) to (2);
  \draw[arc] (1) to (13);
  \draw[arc] (7) to (6);
  \draw[arc] (7) to (9);
  \draw[arc] (8) to (15);
  \draw[arc] (2) to (6);
  \draw[arc] (5) to (8);
  \draw[arc] (5) to (2);
  \draw[arc] (6) to (9);
  \draw[arc] (9) to (13);
  \draw[arc] (10) to (8);
  \draw[arc] (10) to (5);
  \draw[arc] (13) to (10);
\end{tikzpicture}
    \caption{\(H \sqcup K_{3,3}\)}
    \label{fig:order-16-sparse-disconnected-DDMOG}
\end{figure}
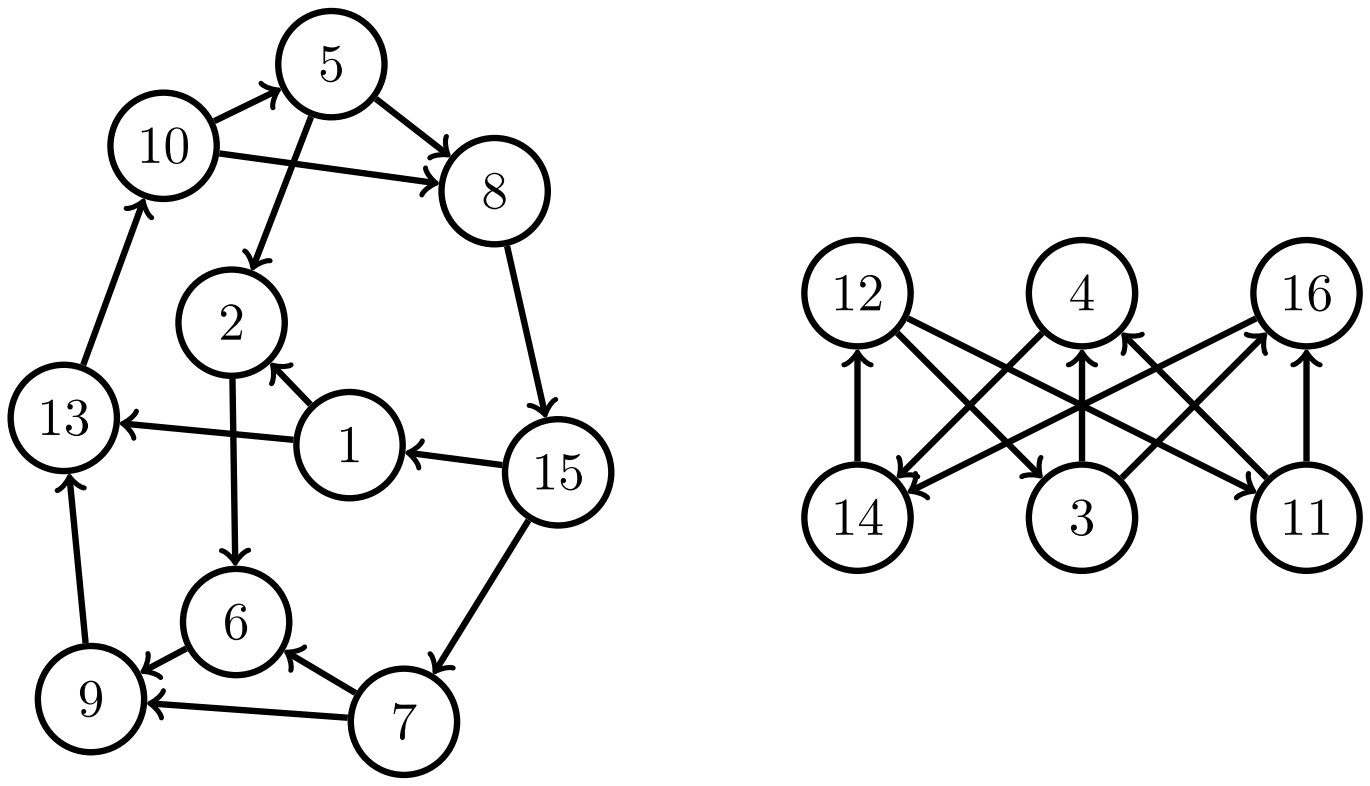

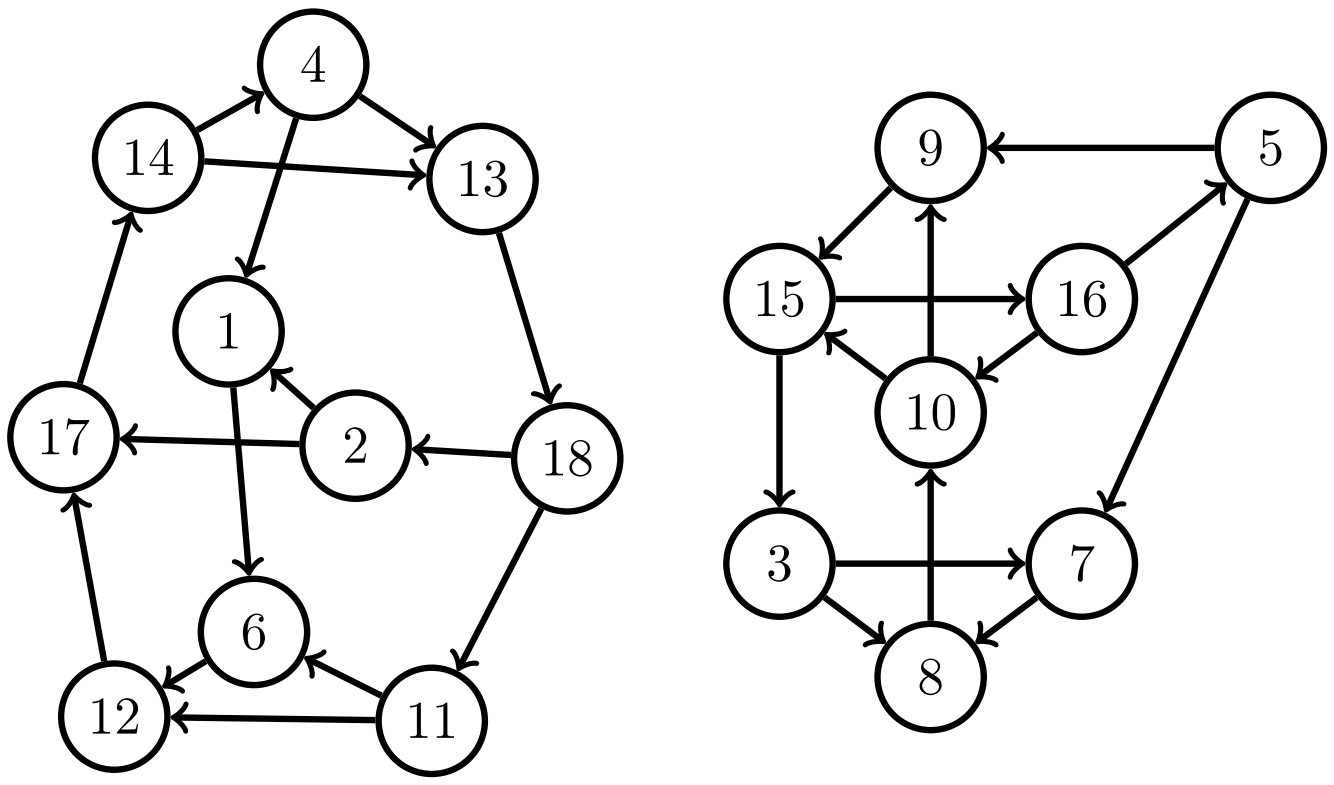
\begin{figure}
    \centering
    \begin{tikzpicture}
    \tikzset{
        labeled_vertex/.style = {circle, draw=black, fill=white, inner sep = 1.5pt, text=black, minimum size=20pt, very thick},
        arc/.style = {->, very thick}
    };
    \begin{scope}[shift={(0,0)},scale=2.8]
        \node[labeled_vertex] (h1) at (0.07, 0.10) {2};
        \node[labeled_vertex] (h3) at (0.25, -0.55) {11};
        \node[labeled_vertex] (h4) at (0.37, 0.73) {13};
        \node[labeled_vertex] (h5) at (-0.23, 0.37) {1};
        \node[labeled_vertex] (h6) at (-0.03, 1.00) {4};
        \node[labeled_vertex] (h7) at (-0.17, -0.34) {6};
        \node[labeled_vertex] (h8) at (-0.50, -0.54) {12};
        \node[labeled_vertex] (h9) at (-0.42, 0.78) {14};
        \node[labeled_vertex] (h10) at (-0.62, 0.12) {17};
        \node[labeled_vertex] (h15) at (0.57, 0.07) {18};
        \draw[arc] (h1) to (h5);
        \draw[arc] (h1) to (h10);
        \draw[arc] (h3) to (h7);
        \draw[arc] (h3) to (h8);
        \draw[arc] (h4) to (h15);
        \draw[arc] (h5) to (h7);
        \draw[arc] (h6) to (h4);
        \draw[arc] (h6) to (h5);
        \draw[arc] (h7) to (h8);
        \draw[arc] (h8) to (h10);
        \draw[arc] (h9) to (h4);
        \draw[arc] (h9) to (h6);
        \draw[arc] (h10) to (h9);
        \draw[arc] (h15) to (h1);
        \draw[arc] (h15) to (h3);
    \end{scope}
    \begin{scope}[shift={(4,-.5)}]
        \node [labeled_vertex] (q11) at (-1, 0) {3};
        \node [labeled_vertex] (q2) at (1, 0) {7};
        \node [labeled_vertex] (q12) at (0, -.75) {8};
        \node [labeled_vertex] (q14) at (-1, 1.75) {15};
        \node [labeled_vertex] (q8) at (1, 1.75) {16};
        \node [labeled_vertex] (q13) at (0, 1) {10};
        \node [labeled_vertex] (q6) at (0, 2.75) {9};
        \node [labeled_vertex] (q1) at (2.25, 2.75) {5};
        
        \draw[arc] (q11) to (q2);
        \draw[arc] (q2) to (q12);
        \draw[arc] (q11) to (q12);
        \draw[arc] (q14) to (q8);
        \draw[arc] (q8) to (q13);
        \draw[arc] (q13) to (q14);
        \draw[arc] (q14) to (q11);
        \draw[arc] (q12) to (q13);
        \draw[arc] (q6) to (q14);
        \draw[arc] (q13) to (q6);
        \draw[arc] (q1) to (q6);
        \draw[arc] (q8) to (q1);
        \draw[arc] (q1) to (q2);
    \end{scope}
\end{tikzpicture}
    \caption{\(H\sqcup Q\)}
    \label{fig:order-18-sparse-disconnected-DDMOG}
\end{figure}

\begin{figure}
    \centering
    \begin{tikzpicture}
    \tikzset{
        labeled_vertex/.style = {circle, draw=black, fill=white, inner sep = 1.5pt, text=black, minimum size=20pt, very thick},
        arc/.style = {->, very thick}
    };
    \begin{scope}[shift={(0,0)},scale=2.8]
        \node[labeled_vertex] (h1) at (0.07, 0.10) {22};
        \node[labeled_vertex] (h3) at (0.25, -0.55) {1};
        \node[labeled_vertex] (h4) at (0.37, 0.73) {23};
        \node[labeled_vertex] (h5) at (-0.23, 0.37) {2};
        \node[labeled_vertex] (h6) at (-0.03, 1.00) {4};
        \node[labeled_vertex] (h7) at (-0.17, -0.34) {26};
        \node[labeled_vertex] (h8) at (-0.50, -0.54) {3};
        \node[labeled_vertex] (h9) at (-0.42, 0.78) {25};
        \node[labeled_vertex] (h10) at (-0.62, 0.12) {27};
        \node[labeled_vertex] (h15) at (0.57, 0.07) {29};
        \draw[arc] (h1) to (h5);
        \draw[arc] (h1) to (h10);
        \draw[arc] (h3) to (h7);
        \draw[arc] (h3) to (h8);
        \draw[arc] (h4) to (h15);
        \draw[arc] (h5) to (h7);
        \draw[arc] (h6) to (h4);
        \draw[arc] (h6) to (h5);
        \draw[arc] (h7) to (h8);
        \draw[arc] (h8) to (h10);
        \draw[arc] (h9) to (h4);
        \draw[arc] (h9) to (h6);
        \draw[arc] (h10) to (h9);
        \draw[arc] (h15) to (h1);
        \draw[arc] (h15) to (h3);
    \end{scope}

    \begin{scope}[shift={(5,0)},scale=2.8]
        \node[labeled_vertex] (hp1) at (0.07, 0.10) {7};
        \node[labeled_vertex] (hp3) at (0.25, -0.55) {8};
        \node[labeled_vertex] (hp4) at (0.37, 0.73) {15};
        \node[labeled_vertex] (hp5) at (-0.23, 0.37) {6};
        \node[labeled_vertex] (hp6) at (-0.03, 1.00) {9};
        \node[labeled_vertex] (hp7) at (-0.17, -0.34) {16};
        \node[labeled_vertex] (hp8) at (-0.50, -0.54) {14};
        \node[labeled_vertex] (hp9) at (-0.42, 0.78) {21};
        \node[labeled_vertex] (hp10) at (-0.62, 0.12) {24};
        \node[labeled_vertex] (hp15) at (0.57, 0.07) {30};
        \draw[arc] (hp1) to (hp5);
        \draw[arc] (hp1) to (hp10);
        \draw[arc] (hp3) to (hp7);
        \draw[arc] (hp3) to (hp8);
        \draw[arc] (hp4) to (hp15);
        \draw[arc] (hp5) to (hp7);
        \draw[arc] (hp6) to (hp4);
        \draw[arc] (hp6) to (hp5);
        \draw[arc] (hp7) to (hp8);
        \draw[arc] (hp8) to (hp10);
        \draw[arc] (hp9) to (hp4);
        \draw[arc] (hp9) to (hp6);
        \draw[arc] (hp10) to (hp9);
        \draw[arc] (hp15) to (hp1);
        \draw[arc] (hp15) to (hp3);
    \end{scope}

    \begin{scope}[shift={(-1.7,-4.3)}]
        \node [labeled_vertex] (v4) at (0, 0) {28};
        \node [labeled_vertex] (v5) at (1.5, 0) {10};
        \node [labeled_vertex] (v6) at (3, 0) {18};
        \node [labeled_vertex] (v1) at (0, 1.5) {5};
        \node [labeled_vertex] (v2) at (1.5, 1.5) {12};
        \node [labeled_vertex] (v3) at (3, 1.5) {17};
        \draw[arc] (v4) to (v1);
        \draw[arc] (v1) to (v5);
        \draw[arc] (v1) to (v6);
        \draw[arc] (v2) to (v4);
        \draw[arc] (v5) to (v2);
        \draw[arc] (v6) to (v2);
        \draw[arc] (v3) to (v4);
        \draw[arc] (v5) to (v3);
        \draw[arc] (v6) to (v3);
    \end{scope}

    \begin{scope}[shift={(3.3,-4.3)}]
        \node [labeled_vertex] (w4) at (0, 0) {32};
        \node [labeled_vertex] (w5) at (1.5, 0) {13};
        \node [labeled_vertex] (w6) at (3, 0) {19};
        \node [labeled_vertex] (w1) at (0, 1.5) {31};
        \node [labeled_vertex] (w2) at (1.5, 1.5) {11};
        \node [labeled_vertex] (w3) at (3, 1.5) {20};
        \draw[arc] (w4) to (w1);
        \draw[arc] (w1) to (w5);
        \draw[arc] (w1) to (w6);
        \draw[arc] (w2) to (w4);
        \draw[arc] (w5) to (w2);
        \draw[arc] (w6) to (w2);
        \draw[arc] (w3) to (w4);
        \draw[arc] (w5) to (w3);
        \draw[arc] (w6) to (w3);
    \end{scope}

\end{tikzpicture}
    \caption{\(2H\sqcup2K_{3,3}\)}
    \label{fig:order-32-sparse-disconnected-DDMOG}
\end{figure}
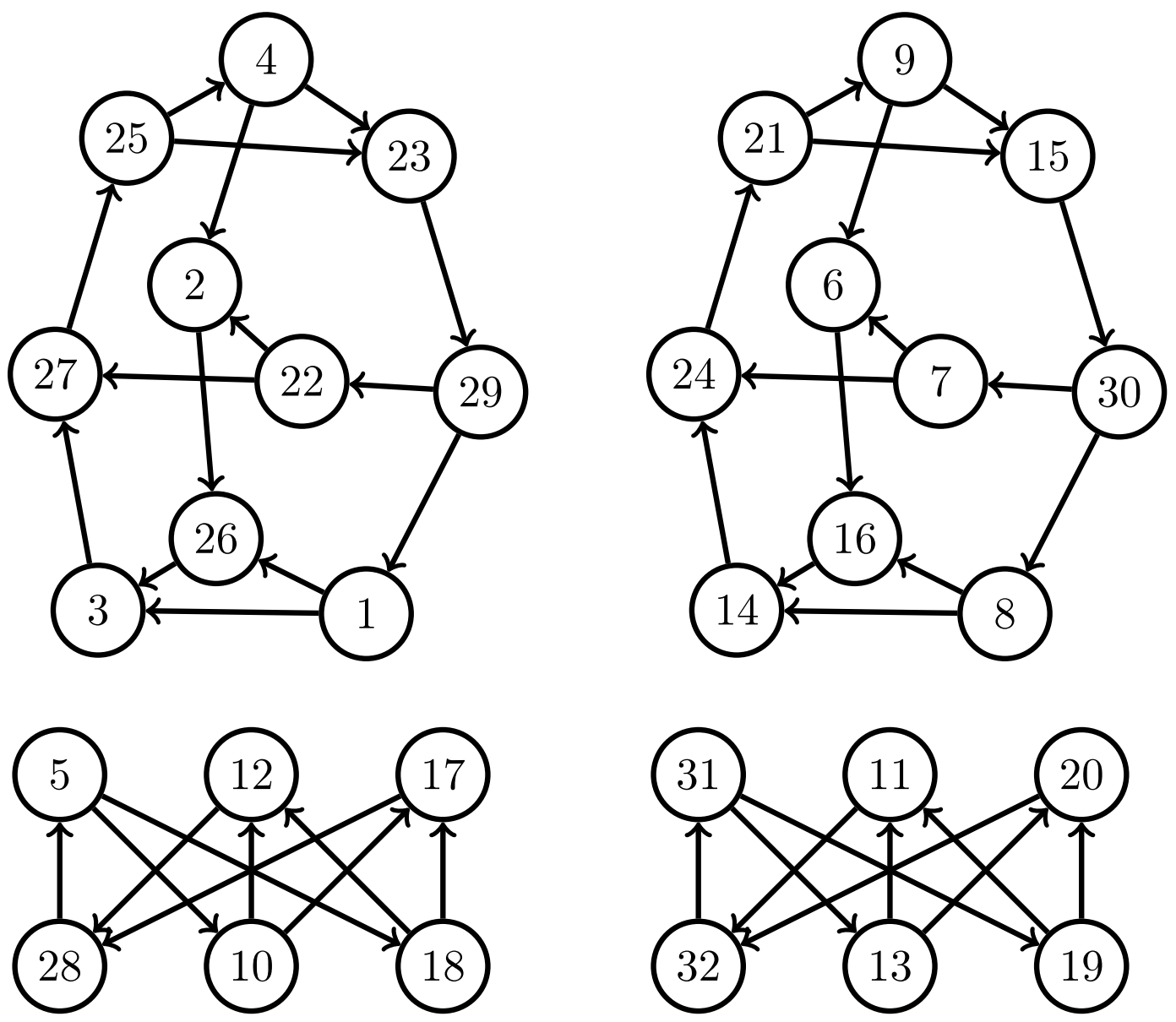

\begin{figure}
    \centering
    \begin{tikzpicture}
    \tikzset{
        labeled_vertex/.style = {circle, draw=black, fill=white, inner sep = 1.5pt, text=black, minimum size=20pt, very thick},
        arc/.style = {->, very thick}
    };
    \node [labeled_vertex] (9) at (0,1) {\(9\)};
    \node [labeled_vertex] (5) at (-1,2) {\(5\)};
    \node [labeled_vertex] (2) at (1,2) {\(2\)};
    \draw[arc] (9) to (5);
    \draw[arc] (5) to (2);
    \draw[arc] (2) to (9);

    \node [labeled_vertex] (3) at (0,0) {\(3\)};
    \node [labeled_vertex] (1) at (-1,-1) {\(1\)};
    \node [labeled_vertex] (10) at (1,-1) {\(10\)};
    \draw[arc] (3) to (1);
    \draw[arc] (10) to (1);
    \draw[arc] (10) to (3);

    \draw[arc] (3) to (9);

    \node [labeled_vertex] (7) at (-2,-0.5) {\(7\)};
    \draw[arc] (1) to (7);
    \draw[arc] (5) to (7);
    
    \node [labeled_vertex] (4) at (2,-0.5) {\(4\)};
    \draw[arc] (10) to (4);
    \draw[arc] (4) to (2);

    \node [labeled_vertex] (6) at (-1.25,-2) {\(6\)};
    \draw[arc] (1) to (6);
    \draw[arc] (7) to (6);

    \node [labeled_vertex] (8) at (1.25,-2) {\(8\)};
    \draw[arc] (8) to (10);
    \draw[arc] (4) to (8);

    \draw[arc] (6) to (8);
\end{tikzpicture}
    \caption{A connected minimal sparsity DDMOG with \(10\) vertices}
    \label{fig:order-10-sparse-connected-DDMOG}
\end{figure}
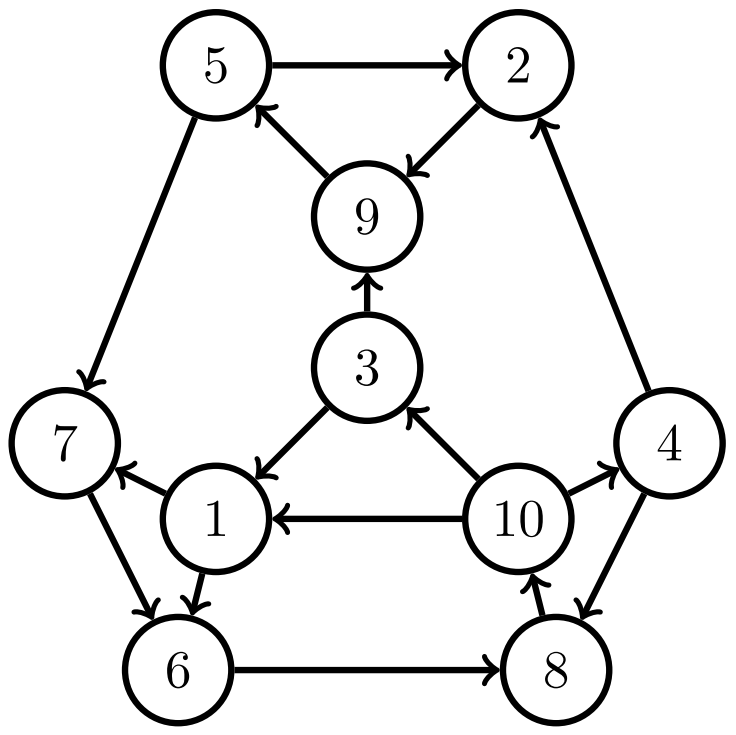

\begin{figure}
    \centering
    \begin{tikzpicture}
    \tikzset{
        labeled_vertex/.style = {circle, draw=black, fill=white, inner sep = 1.5pt, text=black, minimum size=20pt, very thick},
        arc/.style = {->, very thick}
    };
    \node [labeled_vertex] (11) at (-2,0) {\(11\)};
    \node [labeled_vertex] (1) at (0,0) {\(1\)};
    \node [labeled_vertex] (5) at (2,0) {\(5\)};
    \draw[arc] (1) to (11);
    \draw[arc] (5) to (1);
    
    \node [labeled_vertex] (14) at (-1, 1) {\(14\)};
    \node [labeled_vertex] (8) at (1, 1) {\(8\)};
    \draw[arc] (11) to (14);
    \draw[arc] (14) to (1);
    \draw[arc] (1) to (8);
    \draw[arc] (8) to (5);
    
    \node [labeled_vertex] (13) at (-2, 2) {\(13\)};
    \node [labeled_vertex] (9) at (0, 2) {\(9\)};
    \node [labeled_vertex] (7) at (2, 2) {\(7\)};
    \draw[arc] (13) to (11);
    \draw[arc] (9) to (13);
    \draw[arc] (7) to (9);
    \draw[arc] (5) to (7);

    \node [labeled_vertex] (10) at (-1, 3) {\(10\)};
    \node [labeled_vertex] (4) at (1, 3) {\(4\)};
    \draw[arc] (14) to (10);
    \draw[arc] (4) to (8);
    \draw[arc] (4) to (7);

    \node [labeled_vertex] (2) at (-2, 4) {\(2\)};
    \draw[arc] (2) to (13);
    \draw[arc] (10) to (2);

    \node [labeled_vertex] (3) at (-1, 5) {\(3\)};
    \node [labeled_vertex] (12) at (1, 5) {\(12\)};
    \draw[arc] (3) to (2);
    \draw[arc] (3) to (4);
    \draw[arc] (10) to (12);
    \draw[arc] (12) to (4);

    \node [labeled_vertex] (6) at (0, 6) {\(6\)};
    \draw[arc] (6) to (3);
    \draw[arc] (12) to (6);
    \draw[arc] (6) to (9);

\end{tikzpicture}
    \caption{A connected minimal sparsity DDMOG with \(14\) vertices}
    \label{fig:order-14-sparse-connected-DDMOG}
\end{figure}
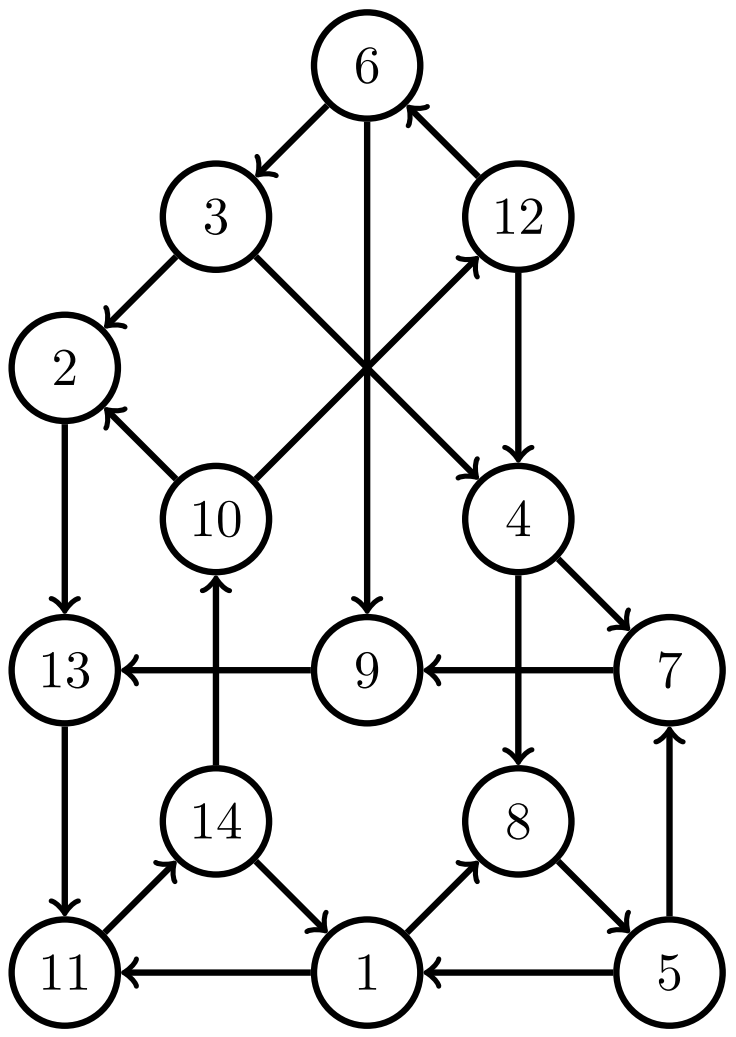

\begin{algorithm}[h!]
\caption{DDMOG Counter}\label{alg:DDMOGIterator}

\begin{algorithmic}[1]
\Require $n$
\Ensure The number of nontrivial DDMOGs with $n$ vertices.
\Statex

\Procedure{CountDDMOGs}{$v, G$}
    \State $count \gets 0$ \label{alg:DDMOGIterator:init-count}
    \If{$v = n$} \label{alg:DDMOGIterator:base-case}
        \If{$G$ is nontrivial}
            \State $count \gets 1$ \label{alg:DDMOGIterator:found-nontrivial-DDMOG}
        \EndIf
    \Else
        \State $w$ $\gets$ weight of $v$ in $G$ \label{alg:DDMOGIterator:compute-weight}
        \For{$tuple$ \textbf{in} \(\{-1,0,1\}^{n-v}\)} \label{alg:DDMOGIterator:begin-arc-search}
            \If{dotproduct($tuple, [v+1, v+2, \cdots, n]$) = $-w$} \label{alg:DDMOGIterator:check-DDM}
    
                \State $G' \gets G$ \label{alg:DDMOGIterator:copy-G}
                \For{$x$ \textbf{in} $tuple$}
                    \State $v'$ $\gets$ \(v\) + index of $x$ in $tuple$
                    \If{$x = 1$}
                    \State Add an arc in $G'$ directed from $v$ to $v'$.
                    \ElsIf{$x = -1$}
                      \State Add an arc in $G'$ directed from $v'$ to $v$.
                    \EndIf
                \EndFor
                \State $count$ $\gets$ $count$ $+$ \Call{countDDMOGs}{$v + 1, G'$} \label{alg:DDMOGIterator:recurse}
            \EndIf
        \EndFor
        \EndIf
    \State \Return $count$ \label{alg:DDMOGIterator:return-count}
\EndProcedure
\Statex
\State \Call{CountDDMOGs}{$1, nK_1$} \label{alg:DDMOGIterator:begin-recursion}

\end{algorithmic}
\end{algorithm}

\end{document}